\documentclass[10pt,twoside]{amsart}
\usepackage{lmodern}
\usepackage[T1]{fontenc}
\usepackage{amsmath,amsfonts,amsthm,amssymb,amscd}
\usepackage{mathrsfs}
\usepackage{latexsym}
\usepackage[pdftex]{graphicx}
\usepackage{enumitem}
\usepackage{bbm}
\usepackage{accents}

\usepackage{geometry}
\newtheorem{corollary}{Corollary}[section]
\newtheorem{lemma}{Lemma}[section]
\newtheorem{theorem}{Theorem}[section]

\usepackage{hyperref}  
\hypersetup{colorlinks=true, linktoc=all,
    linkcolor=blue,
		citecolor=blue,
		urlcolor=cyan,
		bookmarksopen=true
		}
\usepackage{hyperref}  
\newcommand{\lnc}{\mathscr{L}}
\def\B{{\mathscr B}}

\title{Sums of values of non-principal characters over shifted primes}
\author{Zarullo Rakhmonov}
\address{A.Dzhuraev Institute of Mathematics,  National Academy of Sciences of Tajikistan}
\email{zarullo.rakhmonov@gmail.com, zarullo-r@rambler.ru}
\date{}

\begin{document}

\begin{abstract}
For a nonprincipal character $\chi$ modulo $D$, when $x\ge D^{\frac56+\varepsilon}$, $(l,D) = 1$, we prove a nontrivial estimate of the form $\sum_{n\le x}\Lambda (n)\chi (n-l)\ll x\exp\left(-0.6\sqrt{\ln D}\right)$  for the sum of values of $\chi$ over a sequence of shifted primes.
Bibliography: 42 references.
\end{abstract}

\maketitle

\section{Introduction}
I.~M.~Vinogradov's method for estimating exponential sums over primes enabled him to solve a number of arithmetic problems with primes. One of these problems concerns the distribution of values of a nonprincipal character over sequences of shifted primes. In 1938, he proved the following:
\emph{If $q$ is an odd prime, $(l,q)=1$, and  $\chi(a)$ is a nonprincipal character modulo q, then}
$$
T_1(\chi )=\sum_{p\le x}\chi (p-l)\ll x^{1+\varepsilon}\sqrt{\frac{1}{q}+\frac{q}{\sqrt[3]{x}}}.
$$
(see \cite{VinigradovIM-1938}). In 1943, I.~M.~Vinogradov  \cite{VinigradovIM-1943} improved this estimate by showing that
\begin{equation}\label{Formula 0tsenka IMV-1943}
|T_1(\chi )|\ll x^{1+\varepsilon} \left(\sqrt{\frac{1}{q}+\frac{q}{x}} +x^{-\frac{1}{6}}\right).
\end{equation}
For $x\gg q^{1+\varepsilon}$, this estimate is nontrivial and implies an \emph{asymptotic formula for the number of quadratic residues (nonresidues) modulo $q$ of the form $p-l$, $p\le x$}.

A Goldbach number is a positive integer that can be represented as the sum of two odd primes. The problem of the distribution of such numbers
in ``short'' arithmetic progressions has appeared in attempts to solve the binary Goldbach problem. The first result of a conditional nature belongs to Yu.V. Linnik \cite{Linnik-1952}. Under the assumption of the Generalized Riemann Hypothesis, he obtained the upper bound
$$
G(D,l)\leq D \ln^{6}D,
$$
where $G(D,l)$ is the smallest Goldbach number in the arithmetic progression
$$
Dk + l, \qquad k = 0,1,2, \ldots .
$$
This result was further refined by K.Prahar \cite{Prachar_Acta-Arith-29(1976),Prachar_Acta-Arith-44(1984)} and Yu. Wang \cite{WangYuan_SciSinica-1977-20}. Under the same assumptions, they proved that
$$
G(D,l)\leq D (\ln D)^{3+\varepsilon}.
$$
In 1968 M.~Jutila \cite {Jutila} proved an unconditional theorem. Using the estimate (\ref{Formula 0tsenka IMV-1943}), he showed that if $ D $ is an odd prime number, then
$$
G(D,l)\ll D^{\frac{11}{8}+\varepsilon}.
$$

Then I.~M.~Vinogradov obtained a nontrivial estimate for $T_1(\chi )$ in the case when $x\ge q^{0.75+\varepsilon}$, where $q$ is a prime (see \cite{VinigradovIM-1952,VinigradovIM-1953,VinigradovIM-1966}). This result was unexpected. The point is that $T_1(\chi)$ can be expressed as a sum over the zeros of the corresponding Dirichlet $L$-function. Then, assuming that the Generalized Riemann hypothesis is valid for $T_1(\chi )$, one obtains a nontrivial estimate but only for $x~\ge~q^{1+\varepsilon}$.

It seemed that one got what was impossible. In 1973, in this connection Yu.~V.~Linnik \cite {Linnik-1973}  wrote, \emph{``Vinogradov's investigations into the asymptotic behavior of Dirichlet characters are of great importance. As early as 1952, he obtained an estimate for the sum of the Dirichlet characters of shifted primes $T_1(\chi)$, which gave a degree of decrease relative to $x$ as soon as $x>q^{0.75+\varepsilon}$, where $q$ is the modulus of the character. This has fundamental significance, since it is deeper than what can be obtained by direct application of the Generalized Riemann hypothesis, and in this direction seems to carry a deeper truth than that hypothesis (if the hypothesis is valid). Recently A.~A.~Karatsuba has been able to improve this estimate.''}

In 1968, Karatsuba developed a method that allowed him to obtain a nontrivial estimate for short character sums in finite fields of fixed degree \cite{KaratsubaAA-2008,KaratsubaAA-1968,KaratsubaAA-1970-1}. In 1970, he improved this method and, combining it with Vinogradov's method, established the following statement \cite{KaratsubaAA-2008,KaratsubaAA-1970-2,KaratsubaAA-1970-3}: \emph{If $q$ is a prime, $\chi(a)$ is a nonprincipal character modulo $q$, and $x\ge q^{\frac{1}{2}+\varepsilon}$, then}
$$
T_1(\chi )\ll xq^{-\frac{1}{1024}\varepsilon^2}.
$$
Karatsuba used these estimates to find asymptotic formulas for the number of quadratic residues and nonresidues of the form $p+k$ and for the number of products of the form $p(p'+k)$ in an arithmetic progression with increasing difference \cite{KaratsubaAA-1970-4} (see also \cite{KaratsubaAA-2008,KaratsubaAA-1971,KaratsubaAA-1975-1,KaratsubaAA-1975-2,KaratsubaAA-1976,KaratsubaAA-1978}).

The present author generalized estimate (\ref{Formula 0tsenka IMV-1943}) to the case of a composite modulus and proved the following statement \cite{RakhmonovZKh-1986-UMN,RakhmonovZKh-1986-DANRT,RakhmonovZKh-1994-TrMIRAN}: \emph{Let $D$ be a sufficiently large positive integer, $\chi$  be a nonprincipal character modulo $D$, $\chi_q$ be the primitive character generated by $\chi$, and $q_1$ be the product of primes that divide $D$ but do not divide $q$; then}
$$
T_1(\chi )\le x\ln^5x \left(\sqrt{\frac{1}{q}+\frac{q}{x}\tau^2(q_1)} +x^{-\frac{1}{6}}\tau (q_1)\right)\tau(q).
$$
Applying this estimate, he proved that the estimate
$$
G(D,I)\ll D^{c+\varepsilon},
$$
holds for a sufficiently large odd integer $D$ (see \cite{RakhmonovZKh-1986-UMN, RakhmonovZKh-1986-IzvANRT}), where $\varepsilon$  is arbitrarily small positive constant, $c$ is the lower bound of such $a$'s, for which the inequality
$$
\sum_{\chi\,mod\,D}N(\alpha,T,\chi)\ll \left(DT\right)^{2a(1-\alpha )} \left(\ln DT\right)^{A}.
$$
holds for some constant $A>2$. In fact, Huxley's ``density'' theorem \cite{Huxley-1942-Inv.Math} allows one to obtain the estimate $c\le\frac{6}{5}$ by taking  $A=14$.

In 2010, J.~B.~Friedlander, K.~Gong, I.~E.~Shparlinski \cite{Friedlander-Gong-Shparlinski} showed that for a composite $q$ the sum $T(\chi_q )$ can be estimated nontrivially when the length $x$ of the sum is of smaller order than $q$. They proved the following: \emph{For a primitive character $\chi_q$ and any $\varepsilon >0$, there exists a $\delta >0$ such that the following estimate holds for all $x\ge q^{\frac{8}{9}+\varepsilon}$:}
$$
T_1(\chi_q )\ll xq^{-\delta}. \eqno (FGSh)
$$

In 2013, the present author proved the following theorem (see \cite{RakhmonovZKh-2013-DANRT, RakhmonovZKh-2013-IzvSarUniv, RakhmonovZKh-2014-Ch.sbor-15-2(50)}: \emph{If $q$ is a sufficiently large positive integer, $\chi_q$ is a primitive character modulo $q$, $(l,q)=1$, $\varepsilon$ is an arbitrarily small fixed
positive number, and $x\ge q^{\frac{5}{6}+\varepsilon}$, then}
\begin{align*}
T_1(\chi_q )\ll x\exp\left(-\sqrt{\ln q}\right).
\end{align*}
In 2017 Bryce Kerr in \cite{Kerr-2017} obtained an estimate (FGSh), i.e. for $T_1(\chi_q )$ obtained an estimate that gives a power saving which is nontrivial for $x\ge q^{\frac{5}{6}+o(1)}$.

In 2017, the author (see \cite{RakhmonovZKh-2017-TrMIRAN-299, RakhmonovZKh-2017-DANRT-9}) proved the following theorem: \emph{Let $D$ be a sufficiently large positive integer, $\chi$ be a nonprincipal character modulo $D$, $\chi_q$ be the primitive character modulo $q$ generated by $\chi$, $q$ be a cube-free number, $(l,D)=1$, and $\varepsilon$  be an arbitrarily small positive constant. Then, for $x\ge D^{\frac{1}{2}+\varepsilon}$, we have
\begin{align*}
T_1(\chi)
\ll x\exp\left(-0,6\sqrt{\ln D}\right).
\end{align*}}

As was noted above, non-trivial estimates of the sum $ T_1(\chi)$, where $\chi$ is a non-principal character modulo $D$, where $D$ is a prime number, have been applied in problems of finding the least Goldbach number and the distribution of products of shifted primes in short arithmetic progressions. To solve these problems for composite modulo $D$, along with non-trivial estimates of the sum $T_1(\chi)$ for primitive characters, similar estimates are required for imprimitive characters. Therefore, we obviously must consider a non-trivial estimation problem for the sum $T_1(\chi)$, where $\chi$ is a non-principal character modulo a composite number $D$.

In this paper we obtain a non-trivial estimate for the sum $ T(\chi)$ for all non-principal characters modulo a composite number. Let us state the main result.

{\theorem \label{Teorema ocnovnaya} Suppose that $D$ is a sufficiently large positive integer, $\chi$ is a non-principal character modulo $D$, $\chi_q$ be the primitive character modulo $q$ generated by $\chi$,  $(l,D)=1$, $\varepsilon$ is an arbitrarily small positive constant. Then for $x\ge D^{\frac{5}{6}+\varepsilon}$ and $q>\exp\sqrt{2\ln D}$, we have
\begin{align*}
T(\chi)=\sum_{n\le x}\Lambda (n)\chi (n-l)\ll x\exp\left(-0.6\sqrt{\ln D}\right),
\end{align*}
where the constant under the sign $\ll$ depends only on $\varepsilon $.}

\emph{\textbf{Notations:}} In what follows, we shall always assume that $D$ is a sufficiently large natural number, $x$, $l$ are natural numbers, $(l,D)=1$, $\chi$ is a non-principal character modulo $D$, $\chi_q$ is a primitive character modulo $q$, generated by the character $\chi$, $q_1$ is product of prime numbers dividing $D$ but not dividing $q$, hence $(q,q_1)=1$, $\nu$ is divisor of $q_1$, $\lnc=\ln D$, $\lnc_q = \ln q$, $\delta\le 10^{-4}$ is a fixed positive number, $c$ is a fixed positive  number, possibly not the same at all times, $ \omega(q)$ is the number of distinct prime divisors of $q$, which admits well-known estimate
\begin{equation}\label{formula-omega(q)-otsenka}
\omega (q)\le \frac{c_\omega\lnc_q}{\ln \lnc_q}, 
\end{equation}.

In what follows, we shall frequently use the following well-known lemmas.

{\lemma \label{Resheto} Let $f(n)$ be an arbitrary complex-valued function, $u_1\le~x$, $r\ge 1$,
$$
C_r^k=\frac{r!}{k!(r-k)!},\quad \lambda (n)=\sum_{d|n, \ d\le u_1}\mu(d).
$$
Then the following identity holds:
\begin{align*}
\sum_{n\le x}\hspace{-2pt}\Lambda(n)f(n)=&\sum_{k=1}^r\hspace{-2pt}(-1)^{k-1}C_r^k\hspace{-4pt}\sum_{m_1\le u_1}\hspace{-6pt}\mu (m_1)\cdots\hspace{-20pt}\sum_{\substack{m_k\le u_1\\ m_1\cdots m_kn_1\cdots n_k\le x}}\hspace{-20pt}\mu (m_k) \sum_{n_1}\cdots \sum_{n_k}
                             \ln n_1f(m_1n_1\cdots m_km_k)\\
                            & +(-1)^r \sum_{n_1>u_1}\lambda (n_1)\cdots \sum_{\substack{n_r>u_1\\n_1\cdots n_rm\le x}}\lambda (n_r) \sum_{m}
                            \Lambda (m)f(n_1\cdots n_rm).
\end{align*}}
{\sc Proof.} Proved in \cite{RakhmonovZKh-1993-IzvRAN} using Heath-Brown identity (see \cite{Heath-Brown-1982}).

{\lemma \label{Lemma AIVinogradova} Let $F(x,z,b)$ be the number of positive integers less than $x$ and relatively prime to $b$, $b\le x$ and having only prime divisors less than $z$, $\ln x\le z\le x^{\frac{1}{e}}$, $\alpha =\ln z/\ln x$; then for some $|\theta|\le 1$, we have
$$
F(x,z,b) \ll x\prod_{p|b}\left(1-\frac{1}{p}\right)
\exp \left(-\frac{1}{\alpha}\left(\ln \frac{1}{\alpha}+\ln \ln \frac{1}{\alpha}  \right)
+\frac{1}{\alpha}+\frac{2\theta}{\alpha \ln \frac{1}{\alpha}}\right).
$$
}
{\sc Proof} \cite{VinogradovAI}.

{\lemma \label{Teorema Burgess kv} Let $r$ be an arbitrary fixed natural number, $Z$ is a natural number, $q$ is a square-free number or $r=2$. Then the following relations holds:
\begin{align*}
\sum_{\lambda =0 }^{q-1} \left|\sum_{z=1 }^Z\chi_q(\lambda +z)\right|^{2r}\ll Z^rq+Z^{2r}q^{\frac{1}{2}+\delta},
\end{align*}
where the constant under the sign $\ll$ depends only on $r$ and $\delta$.
}

{\sc Proof}  \cite{Burgess-Proc.Lon.M.S (3)12-1962}.

{\lemma \label{Teorema Burgess cub} For an arbitrary natural number $Z\le q^\frac16$ the following relations hold:
\begin{align*}
 \sum_{z_1,\ldots ,z_6=1 }^Z\left|\sum_{\lambda =0 }^{q-1}\chi_q\left(\frac{(\lambda +z_1)(\lambda +z_2)(\lambda +z_3)}
{(\lambda +z_4)(\lambda +z_5)(\lambda +z_6)}\right)\right|\ll Z^3q^{1+\delta}.
\end{align*}
} {\sc Proof}  \cite{Burgess-1986}.

{\lemma \label{Kol-vo u <= U,(u,q)=1} The following asymptotic formula holds for all positive integers $q$ and~$U$:
$$
 \left|\sum_{u=1 \atop (u,q)=1}^U1-\frac{\varphi (q)}{q}U\right|\le 2^{\omega (q)}.
$$
}
{\sc Proof}.  We have
\begin{align*}
 \sum_{u=1 \atop (u,q)=1}^U1&= \sum_{u=1}^U\sum_{d|(u,q)}\mu(d)
=\sum_{d|q}\mu (d)\left[\frac{U}{d}\right]=U\sum_{d|q}\frac{\mu (d)}{d}-\sum_{d\backslash q}\mu (d)\left\{\frac{U}{d}\right\}. \end{align*}
Consequently
\begin{align*}
\left| \sum_{u=1 \atop (u,q)=1}^U1-\frac{\varphi (q)}{q}U\right|&= \left|\sum_{d\backslash q}\mu (d)\left\{\frac{U}{d}\right\}\right|
\le \sum_{d|q}\mu^2(d)=\prod_{p|q}(1+\mu^2(p))=2^{\omega (q)}.
\end{align*}

{\lemma \label{Mardjhanashvili}For $x \ge 2$ we have
\begin{align*}
\sum_{n\le x}\tau_r^k(n)\ll x\left(\ln x\right)^{r^k-1}, \qquad k=1,2.
\end{align*}
}
\textsc{Proof}  \cite{Mardjhanashvili}.


\section{Auxiliary lemmas}

{\lemma \label{Lemma summs s bolshimi delit} Let  $q|D$, then
\begin{align*}
&\sum_{\substack{d| q\\ d>\exp(\sqrt{2\lnc})}}\frac{\mu^2 (d)}{d}\ll \exp\left(-0.7\sqrt{\lnc}\right).
\end{align*}
}
{\sc Proof.} After dividing the interval of summation into intervals of the form $M~<~d~\le~2M$, we obtain less than $\lnc_q$ sums $S(M)$ of the form
$$
S(M)=\sum_{\substack{d| q \\ M<d\le 2M}}\frac{\mu^2 (d)}{d}\ll M^{-1}\sum_{\substack{d| q \\ d\le 2M}}1.
$$
Let $q=p^{\alpha_1}_1p^{\alpha_2}_2\ldots p^{\alpha_t}_t$ be the canonical decomposition of $q$ into prime factors and let $r_i$ denote $i$-th prime number. Obviously, there exists such $k$ that
$$
q'=r_1r_2\ldots r_k\le q <r_1r_2\ldots r_k r_{k+1},  \qquad k\ge t.
$$
By the prime number theorem we have
$$
\ln q'=\sum_{i\le k}\ln r_i=\sum_{p\le r_k}\ln p > \frac{r_k}{2},
$$
thus
$$
r_k< 2 \ln q'\le 2\lnc_q.
$$
Let $q''=r^{\alpha_1}_1r^{\alpha_2}_2\ldots r^{\alpha_t}_t$. Due to $q''\le q$  and $r_t\le r_k$,  we have
 \begin{equation*}
 S(M)\ll M^{-1} \sum_{\substack{d| q\\ d\le 2M}}1\le  M^{-1}\sum_{\substack{d| q'' \\ d\le 2M}}1.
 \end{equation*}
Prime factors of numbers $d$,  $d| q''$ satisfy the condition $r_j\le  r_t$. As $r_t<2\lnc_q$, the latter sum does not exceed the number of such positive integers that are less than  $2M$,
and have only prime divisors less than $2\lnc_q$, i.e.
$$
 S(M)\ll   M^{-1}\sum_{\substack{d| q''\\ d\le 2M}}1\le M^{-1}F(2M,2\lnc_q, 1)\le M^{-1}F(2M,2\lnc, 1).
$$
Applying Lemma \ref{Lemma AIVinogradova} with
$$
x=2M,\qquad b=1,\qquad z=2\lnc, \qquad \alpha =\frac{\ln z}{\ln x}=\frac{\ln2\lnc}{\ln 2M},
$$
we have
\begin{align*}
S&(M)\ll
\exp\left(-\frac{\ln 2M}{\ln2\lnc}\left(\ln \ln 2M-\B\right)\right), \\
& \B=\ln \ln2\lnc-\ln \ln \frac{\ln 2M}{\ln2\lnc} +1+2\theta \left(\ln \frac{\ln 2M}{\ln2\lnc}\right)^{-1}.
\end{align*}
The condition $2M>\exp(\sqrt{2\lnc})$ implies that $2\lnc<(\ln 2M)^2$. Therefore,
\begin{align*}
\B&\le  \ln \ln \left(\ln 2M\right)^2-\ln \ln \frac{\ln 2M}{\ln \left(\ln 2M\right)^2}+1+2\left(\ln \frac{\ln 2M}{\ln \left(\ln 2M\right)^2}\right)^{-1}\\
&=1+ \ln 2-\ln \left(1 -\frac{\ln 2+\ln\ln\ln 2M}{\ln\ln 2M}\right)+\frac{2}{\ln\ln 2M }\left(1 -\frac{\ln 2+\ln\ln\ln 2M}{\ln\ln 2M}\right)^{-1}\\
&=1+ \ln 2+O\left(\frac{\ln 2+\ln\ln\ln 2M}{\ln\ln 2M}\right)<3< 0.002\ln\ln 2M.
\end{align*}
Consequently,
\begin{align*}
S(M)&\ll \exp\left(-\frac{\ln 2M}{\ln2\lnc}\left(\ln \ln 2M-\B\right)\right)\ll \exp\left(-0.998\cdot\frac{\ln 2M\ln \ln 2M}{\ln2\lnc}\right).
\end{align*}
The condition $2M>\exp(\sqrt{2\lnc})$ implies that  $\ln 2M > \sqrt{2\lnc}$, $\ln \ln 2M >0.5\ln 2\lnc$, therefore
\begin{align*}
S(M)&\ll \exp\left(-0.499\sqrt{2}\sqrt{\lnc}\right)\ll \lnc^{-1}\cdot \exp\left(-0.7\sqrt{\lnc}\right).
\end{align*}
This completes the proof of the Lemma.

{\lemma \label{Chislo resheniy sravn} Let $K$ denote the number of solutions of the congruence
\begin{align*}
&(nd +\eta k)y\equiv(n_1 d +\eta k)y_1\hspace{-7pt} \pmod q, \\
& M<n,n_1 \le M+N, \quad 1\le y,y_1\le Y, \quad  (y,q)=1, \quad  (y_1,q)=1,
\end{align*}
where $(\eta ,q)=(k,d)=1$, $d$ divides $q$,   $2NY<q$, $d<Y$ and
$\rho (qd^{-1},Y)$ denotes the number of divisors $\beta$ of the integer $qd^{-1}$, satisfying the conditions
$qY^{-1}\le \beta <qd^{-1}$  and $(\beta,d)=1$.
Then the following relation holds:
\begin{align*}
 K \le   NY + \frac{2Y^2}{d}+\frac{2Y^2}{d}\rho (qd^{-1},Y)
 +\frac{2(NY)^{1+\delta }}{d},
\end{align*}
where $\delta$ is an arbitrarily small positive constant.}

{\sl \textsc{Proof}.} Dividing both sides of the congruence by $y$ when $y=y_1$, we have
\begin{align*}
nd +\eta k\equiv n_1  d +\eta k \hspace{-7pt} \pmod q, \quad M<n,n_1  \le M+N, \quad 1\le y\le Y, \quad  (y,q)=1,
\end{align*}
or
\begin{align*}
nd \equiv n_1  d  \hspace{-7pt} \pmod q, \quad M<n,n_1  \le M+N, \quad 1\le y\le Y, \quad  (y,q)=1.
\end{align*}
Let us divide both sides of the congruence and modulo by $d$, we have
\begin{align*}
n-n_1 \equiv 0  \hspace{-7pt} \pmod{qd^{-1}}, \quad M<n,n_1  \le M+N, \quad 1\le y\le Y, \quad  (y,q)=1.
\end{align*}
It follows from the conditions $|n-n_1|<N$  and $2N\le qY^{-1}< qd^{-1}$ that the latter congruence becomes an equality
\begin{align*}
n- n_1=0. \quad M<n,n_1  \le M+N, \quad 1\le y\le Y, \quad  (y,q)=1,
\end{align*}
i.e. $n=n_1$ if $y=y_1$.
Thus, we obtain
\begin{equation}\label{K = NY_q+2 kappa}
K\le NY +2\kappa ,
\end{equation}
where $\kappa$ is a number of solutions of the congruence
\begin{align*}
&(nd +\eta k)y\equiv(n_1 d +\eta k)y_1\hspace{-7pt} \pmod q, \\
& M<n,n_1 \le M+N, \quad 1\le y<y_1\le Y, \quad  (y,q)=1, \quad  (y_1,q)=1.
\end{align*}
 or the congruence
\begin{align}  \label{sravn v lemm-1}
&(ny -n_1 y_1)d\equiv \eta k(y_1-y)\hspace{-7pt} \pmod q, \\
& M<n,n_1 \le M+N, \quad 1\le y<y_1\le Y, \quad  (y,q)=1, \quad  (y_1,q)=1. \nonumber
\end{align}
Left-hand side of congruence (\ref{sravn v lemm-1}) and it's modulo is divisibly by $d$.
Therefore, its right-hand side, that is, $\eta k(y_1-y)$, is also divisible by $d$.
As $\eta k$ is coprime with $d$, the number $y_1-y$ must divide $d$, i.e. $y_1-y\equiv 0\hspace{-4pt} \pmod d$, or equivalently, $y_1=y+td$. Therefore, the congruence (\ref{sravn v lemm-1}) can be expressed in the form
\begin{align*}
&(n y-n_1(y+td))d\equiv \eta ktd\hspace{-7pt} \pmod q, \\
&  M<n,n_1 \le M+N, \quad 1\le y<y+td\le Y, \quad  (y,q)=1, \quad  (y+td,q)=1.
\end{align*}
Dividing both sides of the congruence and its modulo by $d$, we obtain
\begin{align}  \label{sravn v lemm-2}
&(n -n_1)y\equiv (n_1d+\eta k)t\hspace{-7pt} \pmod{qd^{-1}}, \\
& M<n,n_1 \le M+N, \quad 1\le y<y+td\le Y, \quad  (y,q)=1, \quad  (y+td,q)=1.\nonumber
\end{align}
Let us split the set of solutions of (\ref{sravn v lemm-2})
\begin{equation}\label{kappa=kappa_1+kappa_2+kappa_3l}
    \kappa =\kappa_1+\kappa_2+\kappa_3,
\end{equation}
into three subsets, where $\kappa_1$, $\kappa_2$ and  $\kappa_3$  denote the number of solutions of the congruence  (\ref{sravn v lemm-2})  having, respectively, the following properties:
\begin{enumerate}
  \item $n_1d+\eta k \equiv 0\hspace{-3pt} \pmod{qd^{-1}}$;
  \item $(n_1d+\eta k)t\equiv 0\!\pmod{qd^{-1}}$  and $n_1d+\eta k\not \equiv 0\hspace{-3pt} \pmod{qd^{-1}}$;
  \item $(n_1d+\eta k)t\not \equiv  0\hspace{-3pt} \pmod{qd^{-1}}$.
\end{enumerate}

\textbf{Estimate of $\kappa_1$.}  The congruence $n_1d+\eta k \equiv 0\hspace{-3pt} \pmod{qd^{-1}}$ has no solution for
 $(d,qd^{-1})>1$, and for  $(d,qd^{-1})=1$ has no more than one solution
$n_1=n_1^*$, $(n_1^*,qd^{-1})=1$ as $2N<qd^{-1}$, i.e. $N$ --- the length of the interval of possible values for $n_1$ --- is less than the modulus of congruence. For $n_1=n_1^*$ the congruence (\ref{sravn v lemm-2}) becomes
\begin{align*}
&(n -n_1^*)y\equiv 0\hspace{-7pt} \pmod{qd^{-1}}, \quad M<n\le M+N, \\
&1\le y<y+td\le Y, \quad  (y,q)=1, \quad  (y+td,q)=1,
\end{align*}
and for fixed $y$ and $t$ has a single solution $n=n_1^*$, therefore
$$
\kappa_1 \le Y\left(\frac{Y}{d}+1\right)\le \frac{2Y^2}{d}.
$$

\textbf{Estimate of $\kappa_2$.} Using condition of the case 2, we represent the congruence (\ref{sravn v lemm-2}) as a system of congruences
\begin{align}  \label{sravn v lemm-3}
&(n_1 -n)y\equiv 0\hspace{-7pt} \pmod{qd^{-1}},\\
& (n_1d+\eta k)t\equiv 0\hspace{-3pt} \pmod{qd^{-1}},\nonumber
\end{align}
satisfying conditions
\begin{align*}
 &n_1d+\eta k \not\equiv 0\hspace{-3pt} \pmod{qd^{-1}}, \quad M<n,n_1 \le M+N, \\
 &1\le y<y+td\le Y, \quad  (y,q)=(y+td,q)=1.
\end{align*}
It follows from the conditions $(y,q)=1$, $|n-n_1|<N$  and $2N\le qd^{-1}$ that the first congruence in (\ref{sravn v lemm-3}) is equivalent to the equality $n_1=n$, therefore the number of solutions of the system (\ref{sravn v lemm-3}) equals the number of solutions of the system
\begin{align}  \label{sravn v lemm-4}
&(nd+\eta k)t\equiv 0\hspace{-6pt} \pmod{qd^{-1}},  \quad  nd+\eta k\not\equiv 0\hspace{-3pt} \pmod{qd^{-1}}, \\
&M<n\le M+N, \quad 1\le y<y+td\le Y, \quad  (y,q)=1, \quad  (y+td,q)=1.\nonumber
\end{align}
Product of the numbers $nd+\eta k$ and $t$ is divisible by $qd^{-1}$, but $nd+\eta k$ is not divisible by $qd^{-1}$, therefore
for every solution of the congruence (\ref{sravn v lemm-4}) there exists the divisor $\beta$ of the number $ad^{-1}$, such that $\beta<qd^{-1}$ and
$$
nd+\eta k\equiv 0\hspace{-3pt} \pmod{\beta}, \qquad t\equiv 0\hspace{-3pt} \pmod{q(d\beta)^{-1}}.
$$
It follows from the condition $(d,\eta k)=1$ that the congruence $nd+\eta k \equiv 0\hspace{-3pt} \pmod{\beta}$
has a solution only for $(\beta,d)=1$. For $(\beta,d)=1$ let us denote by $\kappa_2(\beta)$ the number of solutions of the system of congruences
\begin{align*}
&n\equiv \eta d_\beta^{-1}\hspace{-3pt} \pmod{\beta}, \qquad M<n\le M+N,\quad dd_\beta^{-1}\equiv 1\hspace{-3pt} \pmod{\beta}\\
&t\equiv 0\hspace{-3pt} \pmod{q(d\beta)^{-1}}, \quad 1\le y<y+td\le Y, \quad  (y,q)=1, \quad  (y+td,q)=1,
\end{align*}
or the number of solutions of the congruence
\begin{align*}
n\equiv &\eta d_\beta^{-1}\hspace{-3pt} \pmod{\beta}, \qquad M<n\le M+N,\quad dd_\beta^{-1}\equiv 1\hspace{-3pt} \pmod{\beta}\\
 & 1\le y<y+tq\beta^{-1}\le Y, \quad  (y,q)=1, \quad  (y+tq\beta^{-1},q)=1.
\end{align*}
Limits of the variables $y$ and $t$ in the last congruence can be expressed as
\begin{align}
 1\le y<Y, \quad 1\le t\le \frac{Y-y}{q\beta^{-1}}, \quad  (y,q)=1, \quad  (y+tq\beta^{-1},q)=1. \label{granitsa izm y,t 1}
\end{align}
When $y>Y-q\beta^{-1}$, the upper bound of variable $t$ is smaller than the lower bound, therefore the region (\ref{granitsa izm y,t 1}) can be represented as
\begin{align*}
 1\le y\le Y-q\beta^{-1}, \quad 1\le t\le \frac{Y-y}{q\beta^{-1}}, \quad  (y,q)=1, \quad  (y+tq\beta^{-1},q)=1.
\end{align*}
In its turn, if $\beta \le qY^{-1}$, then $Y-q\beta^{-1}$~---~the upper bound of $y$ is smaller than the the lower bound. Therefore,
$\kappa_2(\beta)=0$ when $\beta \le qY^{-1}$, while for $qY^{-1}< \beta <qd^{-1}$ we have the following estimate for $\kappa_2(\beta)$:
\begin{align*}
\kappa_2(\beta)& \le \left( \frac{N}{\beta}+1\right)\sum_{1\le y<Y-q\beta^{-1} \atop (y,q)=1}\left[\frac{Y-y}{q\beta^{-1}}\right]
\le \left( \frac{NY}{q}+\frac{Y\beta}{q}\right)\sum_{1\le y<Y-q\beta^{-1} \atop (y,q)=1}1.
\end{align*}

And using the relations $2NY<q$, $\beta <qd^{-1}$ and $d<Y$, we find that
\begin{align*}
\kappa_2(\beta)& \le \left(1+\frac{Y}{d}\right)\sum_{1\le y<Y-q\beta^{-1} \atop (y,q)=1}1
\le \frac{2Y}{d}\sum_{1\le y<Y-q\beta^{-1} \atop (y,q)=1}1<\frac{2Y^2}{d}.
\end{align*}
Summing the last inequality over divisors  $\beta$ of the number $qd^{-1}$, satisfying the conditions
$qY^{-1}\le \beta <qd^{-1}$  and $(\beta,d)=1$ and denoting the number of such divisors by $\rho (qd^{-1},Y)$, we obtain
\begin{align*}
\kappa_2&\le \sum_{(\beta,d)=1, \ \beta | qd^{-1}\atop qY^{-1}\le \beta <qd^{-1}}\kappa_2 (\beta)\le
\frac{2Y^2}{d}\rho (qd^{-1},Y).
\end{align*}

\textbf{Estimate of $\kappa_3$.} Let us recall that $\kappa_3$ denotes the number of solutions of the congruence
\begin{align*}
&(n_1 -n)y\equiv (n_1d+\eta k)t\hspace{-7pt} \pmod{qd^{-1}},
\end{align*}
satisfying the conditions
\begin{align*}
& (n_1d+\eta k )t\not \equiv  0\hspace{-3pt} \pmod{qd^{-1}}, \quad M<n,n_1 \le M+N, \\
&1\le y<y+td\le Y, \quad  (y,q)=1, \quad  (y+td,q)=1.
\end{align*}
Fox a fixed pair $(n_1^*,t^*)$ we denote by $\kappa_3(\lambda)$  the number of solutions of the congruence
\begin{align}
\label{sravn v lemm-6}
&(n -n_1^* )y\equiv\lambda\hspace{-12pt}\pmod{qd^{-1}}, \quad M<n \le M+N, \quad 1\le y<y+t^*d\le Y, \quad  (y,q)=1,
\end{align}
where  $0<|\lambda|\le q/2d$ is the absolute least residue of the number $(n_1^*d+\eta k)t^*$ modulo $qd^{-1}$. Using the bounds of variables $n$, $n_1$  and $y$ and condition $2NY<q$, we find that
\begin{equation*}
0<|(n -n_1^* )y|<NY< \frac{q}{2}.
\end{equation*}
It follows that the inequality (\ref{sravn v lemm-6}) becomes an equality
\begin{align}
\label{uravn-7}
&(n -n_1^* )y= \lambda , \quad M<n \le M+N, \quad 1\le y<y+t^*d\le Y, \quad  (y,q)=1,
\end{align}
where the parameter $\lambda$  satisfies
$$
1\le |\lambda | < NY.
$$
Thus, for a fixed pair $(n_1^*,t^*)$ the number of solutions $\kappa_3(\lambda)$ of the congruence
(\ref{sravn v lemm-6}) equals the number of solutions of the equation (\ref{uravn-7}),
and satisfies the following inequality
$$
\kappa_3(\lambda)\le\tau (|\lambda|)\le 0.5\left( NY\right)^\delta .
$$
The number of possible pairs $(n_1^*,t^*)$ does not exceed $N\left(Yd^{-1}+1\right)$. Consequently,
$$
\kappa_3\le N\left(\frac{Y}{d}+1\right)\cdot 0.5(NY)^{\delta }\le \frac{(NY)^{1+\delta }}{d}.
$$
Inserting the estimates for $\kappa_1$,  $\kappa_2$  and  $\kappa_3$  into (\ref{kappa=kappa_1+kappa_2+kappa_3l}),
and then into (\ref{K = NY_q+2 kappa}), we find that
\begin{align*}
 K\le&NY+2(\kappa_1+\kappa_2+\kappa_3)\le  NY+\frac{2Y^2}{d}+\frac{2Y^2}{d}\rho (qd^{-1},Y)
 +\frac{2(NY)^{1+\delta }}{d}.
\end{align*}
This completes the proof of the Lemma.


\section{ Estimates of short character sums}
{\lemma \label{Lemma korotkaya summa Ots Bergess} Let $M$, $N$, $d$, $k$  and $\eta$ be natural numbers,  $(\eta ,q)=(d,k)=1$,  and
$$
S=\sum_{M-N<n\le M}\chi_q(nd+\eta k ).
$$
Then for $N<q^{\frac{7}{12}}d^{-\frac{1}{2}}$, $D^{\frac{\scriptstyle 1}{\scriptstyle 2}}\le q\le D$ and $d\le \exp(\sqrt{2\lnc})$, the following estimate holds:
\begin{equation}
\label{formula korotk summa Ots Bergess induction}
|S|\le N^{\frac{2}{3}}q^{\frac{1}{9}+\frac{\delta}{2} }d^{\frac{2}{3}}.
\end{equation}
}

\textsc{Proof.} The inequality (\ref{formula korotk summa Ots Bergess induction}) for the sum  $S$  will be proved by induction on $N$.
For $N\le q^{\frac{1}{3}}$  and $d\le \exp(\sqrt{2\lnc})$ the right side of (\ref{formula korotk summa Ots Bergess induction}) admits an estimate
$$
N^{\frac{2}{3}}q^{\frac{1}{9}+\frac{\delta }{2} }d^{\frac{2}{3}}\ge
N^{\frac{2}{3}}q^{\frac{1}{9}+\frac{\delta}{2} }>
N^{\frac{2}{3}}q^{\frac{1}{9}}\ge
N^{\frac{2}{3}}\left(N^3\right)^{\frac{1}{9} }=N,
$$
i.e. the inequality (\ref{formula korotk summa Ots Bergess induction}) is trivial and we take it as an induction basis.

Further, we assume that
$$
N>q^{\frac{1}{3}}, \qquad  d\le \exp(\sqrt{2\lnc}).
$$
Shifting the interval of summation by $h$, $1\le h\le H<N$, in the sum $S$,   we obtain
\begin{align*}
S&=\sum_{M-N<n\le M}\hspace{-8pt}\chi_q((n+h)d+\eta k )+\sum_{M-N<n\le M-N+h}\hspace{-8pt}\chi_q(nd+\eta k )-\sum_{M<n\le M+h}\hspace{-8pt}\chi_q(nd+\eta k).
\end{align*}
Estimating last two sums using induction hypothesis, we have
\begin{equation*}
 |S|\le\left|\sum_{M-N<n\le M}\chi_q((n+h)d+\eta k)\right|+2H^{\frac{2}{3}}q^{\frac{1}{9}+\frac{\delta }2 }d^{\frac{2}{3}},
\end{equation*}
Setting $h=yz$ in this inequality and summing it over $y$ and  $z$ satisfying
$$
1\le y \le Y, \quad (y,q)=1, \quad 1\le z \le Z,  \quad Y=\left[0.5Nq^{-\frac{1}{6}}d\right], \quad  Z=\left[0.5 q^{\frac{1}{6}}d^{-1}\right],
$$
we obtain
\begin{align*}
|S|&\le(Y_qZ)^{-1}\left|\sum_{1\le y\le Y \atop(y,q)=1}\sum_{1\le z\le Z}\sum_{M-N<n\le M}\chi_q((n+yz)d+\eta k)\right|+2(YZ)^{\frac{2}{3}}q^{\frac{1}{9}+\frac{\delta}2 }d^{\frac{2}{3}},
\end{align*}
where $Y_q$ denotes the number of integers $y\in [1,Y]$ coprime to $q$. Determining the number $y^{-1}$  from the congruence $yy^{-1}\equiv 1\pmod q$, we find that
\begin{align*}
|S|
&\le(Y_qZ)^{-1}\sum_{M-N<n\le M}\sum_{1\le y\le Y\atop(y,q)=1}\left|\sum_{1\le z\le Z}\chi_q((nd +\eta k)y^{-1}+zd)\right|+2^{-\frac{1}{3}}N^{\frac{2}{3}}q^{\frac{1}{9}+\frac{\delta}2 }d^{\frac{2}{3}}.
\end{align*}
We denote by $I(\lambda )$ the number of solutions of the congruence
$$
(nd +\eta k)y^{-1}\equiv \lambda \pmod q, \quad M-N<n\le M, \quad 1\le y\le Y, \quad  (y,q)=1,
$$
and obtain
\begin{align}\label{formula sveden ots S k W}
|S|& \le (Y_qZ)^{-1}W+2^{-\frac{1}{3}}N^{\frac{2}{3}}q^{\frac{1}{9}+\frac{\delta }2 }d^{\frac{2}{3}}, \\
&W=\sum_{\lambda =0 }^{q-1}I(\lambda )\left| \sum_{1\le z\le Z} \chi_q(\lambda +zd)\right|. \nonumber
\end{align}
Cubing both side of the equality, using H\"{o}lder inequality and using the inequality
$$
\sum_{\lambda =1 }^qI(\lambda )\le NY_q,
$$
we obtain
\begin{align*}
W^3&
\le \left(N Y_q\right)^2\sum_{\lambda =0 }^{q-1}I(\lambda )\left| \sum_{1\le z\le Z} \chi_q(\lambda +zd)\right|^3.
\end{align*}
By squaring both sides of the last inequality and applying Cauchy inequality, we have
\begin{align*}
W^6&\le \left(N Y_q\right)^4K V, \qquad  K =\sum_{\lambda =0 }^{q-1}I^2(\lambda ),
\qquad V= \sum_{\lambda =0 }^{q-1}\left| \sum_{1\le z\le Z} \chi_q(\lambda +zd)\right|^6.
\end{align*}
Applying Lemma \ref{Teorema Burgess cub}, we obtain
\begin{align*}
V&=\sum_{\lambda=0}^{q-1}\sum_{z_1,\ldots,z_6=1}^Z\chi\left(\frac{(\lambda+z_1d)(\lambda+z_2d)(\lambda+z_3d)}{(\lambda+z_4d)(\lambda+z_5d)(\lambda+z_6d)}\right)\\
&\le\sum_{z_1,\ldots,z_6=1}^Z\left|\sum_{\lambda=0}^{q-1}\chi\left(\frac{(\lambda+z_1d)(\lambda+z_2d)(\lambda+z_3d)}{(\lambda+z_4d)(\lambda+z_5d)(\lambda+z_6d)}\right)\right|
\\
&\le  \sum_{z_1,\ldots ,z_6=1 }^{Zd}\left|\sum_{\lambda =0 }^{q-1}\chi \left(\frac{(\lambda +z_1)(\lambda +z_2)(\lambda +z_3)}
{(\lambda +z_4)(\lambda +z_5)(\lambda +z_6)}\right)\right|\le Z^3d^3q^{1+\delta }.
\end{align*}
Therefore
\begin{align}\label{formula sveden ots W(d) k K}
W^6&\le \left(N Y_q\right)^4Z^3d^3q^{1+\delta }\cdot K.
\end{align}
The sum $K$ equals the number of solutions of the congruence
\begin{align*}
&(nd +\eta k)y^{-1}\equiv(n_1 d +\eta k)y_1^{-1}\hspace{-7pt} \pmod q, \\
&M<n,n_1 \le M+N, \quad 1\le y,y_1\le Y, \quad  (y,q)=1, \quad  (y_1,q)=1,
\end{align*}
or the congruence
\begin{align*}
&(nd +\eta k)y\equiv(n_1 d +\eta k)y_1\hspace{-7pt} \pmod q, \\
& M<n,n_1 \le M+N, \quad 1\le y,y_1\le Y, \quad  (y,q)=1, \quad  (y_1,q)=1.
\end{align*}
All conditions of the Lemma \ref{Chislo resheniy sravn} are satisfied for this congruence:
\begin{align*}
&2NY=2N\left[0.5Nq^{-\frac{1}{6}}d\right]\le N^2q^{-\frac{1}{6}}d<\left(q^{\frac{7}{12}}d^{-\frac{1}{2}}\right)^2q^{-\frac{1}{6}}d=q, \\
&\frac{Y}{d}=\frac{\left[0.5Nq^{-\frac{1}{6}}d\right]}{d}>\frac{\left[0.5q^{\frac{1}{6}}d\right]}{d}>0.3q^{\frac{1}{6}}>1.
\end{align*}
According to this Lemma, we have
\begin{align*}
 K &\le    NY + \frac{2Y^2}{d}+\frac{2Y^2}{d}\rho (qd^{-1},Y)
 +\frac{2(NY)^{1+\delta }}{d}  \\
&=\frac{2(NY)^{1+\delta }}{d}\left(1+\frac{d}{2(NY)^\delta} 
+\frac{Y}{N(NY)^\delta }\left(\rho (qd^{-1},Y)+1\right) \right).
\end{align*}
And taking into account that
\begin{align*}
&d\le \exp(\sqrt{2\lnc})\le
\exp(2\sqrt{\lnc_q})\le  q^\frac{\delta }{4},\quad  \rho (qd^{-1},Y)+1\le  \tau (q)\le q^\delta,
\\&
(NY)^\delta  \ge (0,1\,N^2q^{-\frac{1}{6}}d)^\delta  
>(0,1)^\delta q^{\frac{\delta }{2}},
\end{align*}
we have
\begin{align*}
 K \le &\frac{2(NY)^{1+\delta  }}{d}\left(1+\frac{10^\delta q^\frac{\delta }{4}}{2q^\frac{\delta }{2}} + \frac{10\,q^\delta }{q^{\frac{1}{6}+\frac{\delta }{2}}}\right)
\le \frac{3(NY)^{1+\delta  }}{d}.
\end{align*}
Inserting this estimate into (\ref{formula sveden ots W(d) k K}), and the right hand side of obtained estimate into (\ref{formula sveden ots S k W}), we obtain
\begin{align}
W^6&\le \left(N Y_q\right)^4Z^3d^3q^{1+\delta}\cdot K
\le  3qN^5YY_q^4Z^3d^2\left(qNY\right)^{\delta }, \nonumber \\
|S|&\le\frac{3^{\frac{1}{6}}q^{\frac{1}{6}}N^{\frac{5}{6}}Y^{\frac{1}{6}}d^{\frac{1}{3}}\left(qNY\right)^{\frac{\delta }{6}}}
{Y_q^{\frac{1}{3}}Z^{\frac{1}{2}}}
+2^{-\frac{1}{3}}N^{\frac{2}{3}}q^{\frac{1}{9}+\frac{\delta }{2} }d^{\frac{2}{3}}. \label{formula otsenka S}
\end{align}
Further, using the Lemma \ref{Kol-vo u <= U,(u,q)=1}  and well-known inequalities
$$
\omega (q)\le \frac{c_\omega\lnc_q}{\ln \lnc_q}, \qquad \frac{\varphi (q)}{2q}\le \frac{c_\varphi}{\ln \lnc_q},
$$
where $c_\omega$  and $c_\varphi$ are absolute constants, we find
\begin{align*}
 &\left|Y_q-\frac{\varphi (q)}{q}Y\right|\le q^{ \frac{c_\omega\ln 2}{\ln \lnc_q}}< \frac{\varphi (q)}{2q}q^{\frac{1}{7}}
<  \frac{\varphi (q)}{2q}\left[0.5Nq^{-\frac 16} \right]=\frac{\varphi (q)}{2q}Y,
\end{align*}
i.e.
$$
Y_q>\frac{\varphi (q)}{2q}Y\ge \frac{c_\varphi Y}{\ln \lnc_q}.
$$
Inserting the last  inequality in (\ref{formula otsenka S}), we express the parameter $Y_q$ in terms of $Y$. We have
\begin{align*}
|S|& \le  \frac{3^{\frac{1}{6}}}{c_\varphi^{\frac{1}{6}}}\cdot   \frac{q^{\frac{1}{6}}N^{\frac{5}{6}}d^{\frac{1}{3}}\left(qNY\right)^{\frac{\delta }{6}}(\ln\lnc_q)^{\frac{1}{6}}}{Y^{\frac{1}{6}}Z^{\frac{1}{2}}}
+2^{-\frac{1}{3}}N^{\frac{2}{3}}q^{\frac{1}{9}+\frac{\delta }2 }d^{\frac{2}{3}}.
\end{align*}
Taking into account that $Y=\left[0.5Nq^{-\frac{1}{6}}d\right]$, $Z=\left[0.5 q^{\frac{1}{6}}d^{-1}\right]$ and $N<q^{\frac{7}{12}}d^{-\frac{1}{2}}$,  we find
\begin{align*}
|S|& \le  \frac{3^{\frac{1}{6}}}{c_\varphi^{\frac{1}{6}}}\cdot   \frac{q^{\frac{1}{6}}N^{\frac{5}{6}}d^{\frac{1}{3}}\left(0.5q^{\frac{5}{6}}N^2d\right)^{\frac{\delta }{6}}(\ln\lnc_q)^{\frac{1}{6}}}{(0.25Nq^{-\frac{1}{6}}d)^{\frac{1}{6}}\ (0.25 q^{\frac{1}{6}}d^{-1})^{\frac{1}{2}}}
+2^{-\frac{1}{3}}N^{\frac{2}{3}}q^{\frac{1}{9}+\frac{\delta }2 }d^{\frac{2}{3}}\\
&
=\frac{2^{\frac{2}{3}-\frac{\delta}{6}}\,3^{\frac{1}{6}}}{c_\varphi^{\frac{1}{6}}}\cdot q^{\frac{1}{9}}N^{\frac{2}{3}} d^{\frac{2}{3}} \left(q^{\frac{5}{6}}N^2d\right)^{\frac{\delta }{6}}(\ln\lnc_q)^{\frac{1}{6}}
+2^{-\frac{1}{3}}N^{\frac{2}{3}}q^{\frac{1}{9}+\frac{\delta }2 }d^{\frac{2}{3}} \\
&
\le\frac{2^{\frac{2}{3}-\frac{\delta}{6}}\,3^{\frac{1}{6}}}{c_\varphi^{\frac{1}{6}}}\cdot   N^{\frac{2}{3}}q^{\frac{1}{9}}d^{\frac{2}{3}}q^{\frac{\delta }{3}}(\ln\lnc_q)^{\frac{1}{6}}
+2^{-\frac{1}{3}}N^{\frac{2}{3}}q^{\frac{1}{9}+\frac{\delta }2 }d^{\frac{2}{3}} \\
&
= \left(2^{\frac{2}{3}-\frac{\delta}{6}}\,3^{\frac{1}{6}}c_\varphi^{-\frac{1}{6}}q^{-\frac{\delta }{6}}(\ln\lnc_q)^{\frac{1}{6}}+2^{-\frac{1}{3}}\right)
N^{\frac{2}{3}}q^{\frac{1}{9}+\frac{\delta}2 }d^{\frac{2}{3}}\le N^{\frac{2}{3}}q^{\frac{1}{9}+\frac{\delta }2 }d^{\frac{2}{3}}.
\end{align*}
This completes the proof of the Lemma.

{\lemma \label{Lemma korotkaya summa}  Let $(\eta\nu ,q)=1$,  $y\ge q^{\frac{1}3+\frac{8}{5}\delta}$, $y\le q$, $D^{\frac{\scriptstyle 1}{\scriptstyle 2}}\le q\le D$ and $\nu\le \exp(\sqrt{2\lnc})$, then
\begin{align*}
S_y(u,\eta,\nu)=&\sum_{\substack{u-y<n\le u \\(n,q)=1,\,n\equiv\eta \hspace{-8pt} \pmod\nu}}\chi_q(n-\eta )\ll \frac{ y}{\nu} \exp\left(-0.7\sqrt{\lnc}\right).
\end{align*}
}

{\sc Proof.}   We have an identity
\begin{align*}
S_y(u,\eta,\nu)&=\sum_{d|  q}\mu(d)\sum_{\substack{u-y<nd\le u\\ nd\equiv\eta \hspace{-8pt} \pmod\nu}}\chi_q(nd-\eta).
\end{align*}
Determining the number $d^{-1}_\nu$    from the congruence $dd^{-1}_\nu \equiv1(mod\,\nu)$, we can rewrite the congruence $nd\equiv\eta \hspace{-4pt} \pmod \nu$ as $n\equiv\eta d_\nu^{-1}(mod\,\nu)$. Further, representing the variable of the summation $n$ as $n=\eta d_\nu^{-1} +m\nu$, we have
\begin{align*}
S_y(u,\eta,\nu)&=\sum_{d|  q}\mu(d)\sum_{u-y<(\eta d_\nu^{-1} +m\nu)d\le u}\chi_q((\eta d_\nu^{-1} +m\nu)d-\eta)\\
&= \sum_{d|  q}\mu(d)\sum_{\frac{u-y-\eta dd^{-1}_\nu}{d\nu}<m\le \frac{u-\eta dd^{-1}_\nu}{d\nu}}\chi_q(m\nu d+\eta(dd^{-1}_\nu-1)).
\end{align*}
Representing the congruence $dd^{-1}_\nu \equiv1(mod\,\nu)$ in the form $dd^{-1}_\nu -1=\nu k$, where  $k=k(d,\nu)$ is uniquely defined for each pair of  $d$ and $\nu$, we have
\begin{align}
S_y(u,\eta,\nu)&=\chi_q(\nu)\sum_{d|  q}\mu(d)S_y(u,\eta,\nu,d)\ll \sum_{d|  q}\mu^2(d)|S_y(u,\eta,\nu,d)|,\label{formula Sy(u,eta,nu)<summ}\\
S_y(u,\eta,\nu,d)&=\sum_{u_1-y_1<m\le u_1}\hspace{-14pt}\chi_q(md+\eta k),\quad  u_1=\frac{u-\eta}{d\nu}-\frac{\eta k}{d}, \quad y_1=\frac{y}{d\nu}. \nonumber
\end{align}
Part of the sum in the right side of (\ref{formula Sy(u,eta,nu)<summ}), corresponding to the terms satisfying $d\le\exp(\sqrt{2\lnc})$, will be denoted by $S_y'(u,\eta,\nu)$. Remaining terms are estimated using the trivial estimate for the sum $S_y(u,\eta,\nu,d)$ and the Lemma \ref{Lemma summs s bolshimi delit}:
\begin{align*}
\sum_{\substack{d|  q\\d>\exp(\sqrt{2\lnc})}}\hspace{-12pt}\mu^2(d)&|S_y(u,\eta,\nu,d)|\le \sum_{\substack{d| q\\ \exp(\sqrt{2\lnc})<d}}\hspace{-12pt}\mu^2(d)\left(\frac{y}{\nu d}+1\right)\\
&\le\frac{y}{\nu}\sum_{\substack{d| q\\ \exp(\sqrt{2\lnc})<d}}\hspace{-12pt}\frac{\mu^2(d)}{d}+\tau(q)\ll \frac{y}{\nu}\exp\left(-0.7\sqrt{\lnc}\right).
\end{align*}
Therefore,
\begin{equation}\label{formula S_y(u,eta,nu)=sumS_y(u,eta,nu,d)}
S_y(u,\eta,\nu)\ll\left|S_y'(u,\eta,\nu)\right|+\frac{y}{\nu}\exp\left(-0.7\sqrt{\lnc}\right).
\end{equation}

To estimate $|S_y'(u,\eta,\nu)|$, we consider the following cases:
\begin{enumerate}
  \item $y>\sqrt{q}\,\exp \left( \frac{c_\omega \, \lnc_q}{\ln \lnc_q}\right)$;
  \item $q^{\frac{1}3+\frac{8}{5}\delta}\le y\le \sqrt{q}\,\exp \left( \frac{c_\omega \, \lnc_q}{\ln \lnc_q}\right)$.
\end{enumerate}

\textbf{1.} To $S_y(u,\eta,\nu,d)$ we apply a formula that establishes a relation between the values of primitive
characters and the values of Gauss sums (see \cite{KaratsubaOATCh}, Ch. VIII, \S  1, Lemma 3), which yields
\begin{align*}
S_y(u,\eta,\nu,d)&=\sum_{a=1}^q \sum_{u_1-y_1<m\le u_1,\atop a\equiv md+\eta k \hspace{-7pt}\pmod q}\chi_q(md+\eta k)\\
&=\sum_{a=1}^q \chi_q(a)\sum_{u_1-y_1<m\le u_1}\frac{1}{q}\sum_{t=0}^{q-1}e\left(\frac{(a-md-\eta k )t}{q}\right)\\
&=\frac{\tau (\chi_q)}{q}\sum_{t=0}^{q-1}\bar \chi_q(t)e\left(\frac{-\eta kt}{q}\right)\sum_{u_1-y_1<m\le u_1}e\left(\frac{-mdt}{q}\right).
\end{align*}
It follows from the conditions $\nu\le \exp(\sqrt{2\lnc})$, $d\le \exp(\sqrt{2\lnc})$ and the hypothesis of the present case that the last sum over m is nonempty. Using the equality
\begin{align*}
 \sum_{m_1 \le m\le m_2} e\left(-\frac{mdt }{q}\right)
 = \frac{ \sin \frac{\pi td (m_2-m_1 +1)}{q}}{\sin \frac{\pi td }{q}} e\left(-\frac{td (m_1+m_2)}{2q}\right)
\end{align*}
with integers $m_1$ and $m_2$, passing to estimates, and taking into account that $|\tau (\chi)|=\sqrt{q}$, we find
\begin{align*}
|S_y(u,\eta,\nu,d)|&\le  \frac{1}{\sqrt{q}}\sum_{t=0\atop (t,q)=1 }^{q-1} \left|\sin \frac{\pi t}{q/d}\right|^{-1}\\
&=\frac{1}{\sqrt{q}}\sum_{t_1=0}^{d-1}\sum_{t_2=0\atop (q/dt_1+t_2,q)=1 }^{q/d-1}\left|\sin \frac{\pi t_2}{q/d}\right|^{-1}\le\frac{d}{\sqrt{q}}\sum_{t=1}^{q/d-1} \left|\sin \frac{\pi t}{q/d}\right|^{-1}.
\end{align*}
If $q/d$ is an odd number, then
\begin{align*}
|S_y(u,\eta,\nu,d)|&\le  
\frac{2d}{\sqrt{q}}\sum_{t=1}^{\frac{q/d -1}{2}}\frac{1}{\left|\sin \frac{\pi t}{q/d }\right|}\le
\frac{2d}{\sqrt{q}}\sum_{t=1}^{\frac{q/d -1}{2}}\frac{q/d}{2t }=\sqrt{q}\sum_{t=1}^{\frac{q/d -1}{2}}\frac{1}{t},
\end{align*}
since $\sin\pi\alpha\ge 2\alpha$ for $0\le \alpha\le 1/2$. Using the inequality $\frac{1}{t}\le \ln \frac{2t+1}{2t-1}$, we then find
\begin{align*}
|S_y(u,\eta,\nu,d)|&\le \sqrt{q}\sum_{t=1}^{\frac{q/d -1}{2}}\left(\ln (2t+1)- \ln (2t-1)\right)=\sqrt{q}\ln q/d\le \sqrt{q}\lnc_q.
\end{align*}
If $q/d $ is an even number, then
\begin{align*}
|S_y(u,\eta,\nu,d)|&\le \frac{2d}{\sqrt{q}}\sum_{t=1}^{\frac{q/d }{2}-1}\frac{1}{\left|\sin\frac{\pi t}{q/d}\right|}+\frac{d}{\sqrt{q}}\le\frac{2d}{\sqrt{q}} \sum_{t=1}^{\frac{q/d }{2}-1}\frac{q/d}{2t }+\frac{d}{\sqrt{q}}\ll \sqrt{q}\lnc_q.
\end{align*}
From this and from the definition of $S_y'(u,\eta,\nu)$ we have
\begin{align*}
S_y'(u,\eta,\nu)&\ll\sum_{\substack{d|  q\\d\le \exp(\sqrt{2\lnc})}}\mu^2(d)\sqrt{q}\lnc_q\le\sqrt{q}\lnc_q\sum_{d|  q}\mu^2(d)=2^{\omega(q)}\sqrt{q}\lnc_q.
\end{align*}
Using the inequality (\ref{formula-omega(q)-otsenka}) and relations $\lnc\le 2\lnc_q$, $\nu\le \exp(\sqrt{2\lnc})\le \exp(2\sqrt{\lnc_q})$, we find
\begin{align*}
S_y'(u,\eta,\nu)&\ll\frac{\sqrt{q}\,\nu\,\exp\left(\omega(q)\ln 2+\ln\lnc_q+0.7\sqrt{\lnc}\right)}{y} \cdot \frac{y}{\nu} \exp\left(-0.7\sqrt{\lnc}\right)\\
 &\le \frac{\sqrt{q} \cdot \exp \left(c_\omega\ln 2\ \frac{\lnc_q}{\ln\lnc_q}+3\sqrt{\lnc_q}+\ln\lnc_q\right)}{y}\cdot \frac{y}{\nu} \exp\left(-0.7\sqrt{\lnc}\right)\\
 &\le \frac{\sqrt{q} \cdot \exp \left( \frac{c_\omega \, \lnc_q}{\ln\lnc_q}\right)}{y}\cdot \frac{y}{\nu}\exp\left(-0.7\sqrt{\lnc}\right)\ll \frac{ y}{\nu} \exp\left(-0.7\sqrt{\lnc}\right).
 \end{align*}

\textbf{2.} If $(d,k)>1$, then $S_y(u,\eta,\nu,d)=0$; therefore, without loss of generality, we will assume that $(d,k)=1$ in $S_y(u,\eta,\nu,d)$ and applying Lemma \ref{Lemma korotkaya summa Ots Bergess}, obtain
$$
|S_y(u,\eta,\nu,d)|\le\left(\frac{y}{d\nu}\right)^{\frac{2}{3}}q^{\frac{1}{9}+\frac{\delta }2 }d^{\frac{2}{3}} \le\left(\frac{y}{\nu}\right)^{\frac{2}{3}}q^{\frac{1}{9}+\frac{\delta }2}.
$$
From this and from the definition of $S_y'(u,\eta,\nu)$ we have
\begin{align*}
\left|S_y'(u,\eta,\nu)\right|&\ll\left(\frac{y}{\nu}\right)^{\frac{2}{3}}q^{\frac{1}{9}+\frac{\delta }2}\sum_{\substack{d|  q\\d\le \exp(\sqrt{2\lnc})}}\mu^2(d)\ll\left(\frac{y}{\nu}\right)^{\frac{2}{3}}q^{\frac{1}{9}+\frac{\delta }2}2^{\omega(q)}.
\end{align*}
Using the inequality (\ref{formula-omega(q)-otsenka}) and relation $\lnc\le 2\lnc_q$, we find
\begin{align*}
&\left|S_y'(u,\eta,\nu)\right|\ll\left(\frac{y}{\nu}\right)^{\frac{2}{3}}q^{\frac{1}{9}+\frac{\delta }2}\exp \left(c_\omega\ln 2\ \frac{\lnc_q}{\ln\lnc_q}\right)\\
&
=\frac{ y}{\nu} \exp\left(-0.7\sqrt{\lnc}\right)\cdot \left(\frac{q^{\frac{1}{3}+\frac{3\delta}{2}}\nu\cdot
\exp\left(3c_\omega\ln 2\ \frac{\lnc_q}{\ln\lnc_q}+2,1\sqrt{\lnc}\right)}{y}\right)^{\frac{1}{3}}\\
&
<\frac{ y}{\nu} \exp\left(-0.7\sqrt{\lnc}\right)\cdot \left(\frac{q^{\frac{1}{3}+\frac{3\delta}{2} }\cdot
\exp\left(3c_\omega\ln 2\ \frac{\lnc_q}{\ln\lnc_q}+4\sqrt{\lnc_q}\right)}{y}\right)^{\frac{1}{3}}\\
&< \frac{ y}{\nu} \exp\left(-0.7\sqrt{\lnc}\right)\cdot \left(\frac{q^{\frac{1}{3}+\frac{8\delta}{5}}}{y}\right)^{\frac{1}{3}}\ll \frac{ y}{\nu} \exp\left(-0.7\sqrt{\lnc}\right).
\end{align*}
Inserting our estimate of $|S_y'(u,\eta,\nu)|$ in (\ref{formula S_y(u,eta,nu)=sumS_y(u,eta,nu,d)}) completes the proof of the Lemma.


\section{Estimates of double character sums}
{\lemma \label{Lemma dvoynie summi kv osredn} Let  $M$, $N$, $U$ be integers, $N\le U<2N$, $a_m$  and $b_n$ are such integer-valued functions that
$$
\sum_{M<m\le 2M}|a_m|^\alpha\ll M\lnc^{c_\alpha},\quad  \alpha =1,2; \qquad |b_n|\ll B.
$$
Then the following estimate holds
\begin{align*}
W=\sum_{M<m\le 2M}\hspace{-10pt} a_m\sum_{\substack{U<n\le \min(xm^{-1},2N)\\ (mn,q)=1,\,mn\equiv l\hspace{-8pt} \pmod \nu}} \hspace{-30pt}b_n\chi_q (mn-l)\ll  B\left(M^\frac{3}{4}N^\frac{1}{2}q^{\frac{1}{4}}+M^\frac{3}{4}Nq^{\frac{1}{8}+\frac{\delta}{4}}\right)
\lnc^{\frac{2c_1+c_2}{4}+1}.
\end{align*}}

\textsc{Proof.} It follows from the condition $(l,\nu)=1$ that $(mn,\nu)=1$. Denoting the inner sum in $W$ by $\B(m)$ and representing the congruence $mn\equiv l\hspace{-5pt}\pmod\nu$ in the form $n\equiv lm_\nu^{-1}\hspace{-4pt}\pmod\nu$, we shall transform other sum so that it does not depend on $m$. We have the identity
\begin{align*}
\B(m)&=\sum_{\substack{N<n\le 2N\\ (n,q)=1}}b_n\chi_q(mn-l)\sum_{U<r\le \min(xm^{-1},2N)}\\
& \qquad \qquad \times\frac{1}{q}\sum_{k=0}^{q-1}e\left(\frac{k(n-r)}{q}\right)\frac{1}{\nu}\sum_{j=0}^{\nu -1}e\left(\frac{j(n-lm_\nu^{-1})}{\nu}\right)\\
&=\frac{1}{q\nu}\sum_{j=0}^{\nu -1}e\left(-\frac{lm_\nu^{-1}j}{\nu}\right)\sum_{k=0}^{q-1}\B(k\nu+jq,m)\sum_{U<r\le \min(xm^{-1},2N)}e\left(-\frac{kr}{q}\right),\\
& \qquad \qquad  \B(k\nu+jq,m)=\sum_{\substack{N<n\le 2N\\ (n,q)=1}}b_n\chi_q(mn-l)e\left(\frac{(k\nu+jq)n}{q\nu}\right).
\end{align*}
Denoting $N'=\min([xm^{-1}],2N)$, extracting the term with $k=0$ and summing over $r$, we obtain
\begin{align*}
|\B(m)|\le&\frac{1}{q\nu}\sum_{j=0}^{\nu -1}e\left(-\frac{lm_\nu^{-1}j}{\nu}\right)\left((N'-U)|\B(jq,m)|\right.\\
&\left.+\sum_{k=1}^{q-1}|\B(k\nu+jq,m)|\frac{ \sin \frac{\pi k(N'-U)}{q}}{\sin \frac{\pi k}{q}} e\left(-\frac{k(N'+1+U)}{2q}\right)\right).
\end{align*}
Passing to the inequalities, we find
\begin{align*}
|\B(m)|\le\frac{1}{\nu}&\sum_{j=0}^{\nu -1}\left(\frac{N'-N}{q}|\B(jq,m)|+\frac{1}{q}\sum_{k\le \left[q/2\right]} \frac{|\B(k\nu+jq,m)|}{\left|\sin \frac{\pi k}{q}\right|}\right.\\
&\left.+\frac{1}{q}\sum_{q/2<k\le q-1} \frac{|\B(k\nu+jq,m)|}{\left|\sin \frac{\pi (q-k)}{q}\right|}\right).
\end{align*}
Next, using the condition $N'-N<q$, inequalities $\sin \pi \alpha\ge \alpha$, $0\le \alpha \le 0.5$ and $\frac{1}{k}\le \frac{2}{k+1}$, where $k\ge1$ is an integer, we have
\begin{align*}
&|\B(m)|\le\frac{1}\nu\sum_{j=0}^{\nu -1}\left(\B(jq,m)|+\sum_{k\le\left[q/2\right]}\hspace{-6pt}\frac{|\B(k\nu+jq,m)|}{k}+\sum_{q/2<k\le q-1}\hspace{-10pt}\frac{|\B(k\nu+jq,m)|}{q-k}\right)\\
&\le\frac{2}{\nu}\sum_{j=0}^{\nu -1}\sum_{k=0}^{q-1}\left(\frac{1}{k+1}+\frac{1}{q-k}\right)|\B(k\nu+jq,m)|\ll \lnc \max_{0\le j<\nu}\max_{0\le k<q}|\B(k\nu+jq,m)|.
\end{align*}
From this and from the definition of $W$, we have
\begin{align}
|W|&\ll \lnc \max_{0\le j<\nu}\max_{0\le k<q}W(j,k), \quad W(j,k)=\sum_{\substack{M<m\le 2M\\ (m,q)=1}}|a_m||\B(k\nu+jq,m)|. \label{formula sveden W k W(j,k)}
\end{align}
Let us estimate $W(j,k)$. Squaring both sides of the last equality and applying H\"{o}lder inequality, we have
\begin{align*}
W^2(j,k)&\le\sum_{\substack{M<m\le 2M\\ (m,q)=1}}|a_m|\sum_{\substack{M<m\le 2M\\ (m,q)=1}}\hspace{-6pt}|a_m||\B(kd+jq,m)|^2\\
&\ll M\lnc^{c_1}\sum_{\substack{M<m\le 2M\\ (m,q)=1}}\hspace{-6pt}|a_m||\B(kd+jq,m)|^2.
\end{align*}
Now, squaring both sides of the last inequality and applying Cauchy inequality, we find
\begin{align*}
W^4(j,k)&\ll M^2\lnc^{2c_1}\sum_{\substack{M<m\le 2M\\ (m,q)=1}}|a_m|^2\sum_{\substack{M<m\le 2M\\ (m,q)=1}}|\B(kd+jq,m)|^4\\
&\ll M^3\lnc^{2c_1+c_2}\sum_{\substack{m=0\\(m,q)=1}}^{q-1}\left|\sum_{\substack{N<n\le 2N\\ (n,q)=1}}b_n\chi_q(n-lm_q^{-1})e\left(\frac{(kd+jq)n}{qd}\right)\right|^4\\
&=M^3\lnc^{2c_1+c_2}\sum_{\substack{\lambda=0\\ (\lambda,q)=1}}^{q-1}\left|\sum_{\substack{N<n\le 2N\\ (n,q)=1}}\hspace{-6pt}b_n\chi_q(n+\lambda)e\left(\frac{(kd+jq)n}{qd}\right)\right|^4\\
&\ll M^3\lnc^{2c_1+c_2}\hspace{-12pt}\sum_{\substack{N'<n_1,_2,n_3,n_4\le 2N\\ (n_1,_2,n_3,n_4),q)=1}}\hspace{-8pt}|b_{n_1}b_{n_2}b_{n_3}b_{n_4}|
\left|\sum_{\lambda =0}^{q-1}\chi \left(\frac{(\lambda +n_1)(\lambda +n_2)}{(\lambda +n_3)(\lambda +n_4)}\right)\right|.
\end{align*}
Next, using the estimate  $|b_n|\ll B$, and the Lemma \ref{Teorema Burgess kv}, we obtain
\begin{align*}
&W^4(j,k) \ll M^3\lnc^{2c_1+c_2} B^4\sum_{1\le n_1,\ldots n_4\le 2N}\left|\sum_{\lambda =0}^{q-1}\chi \left(\frac{(\lambda +n_1)(\lambda +n_2)}
{(\lambda +n_3)(\lambda +n_4)}\right)\right| \\
&
\ll M^3\lnc^{2c_1+c_2} B^4\left(N^2q+N^4q^{\frac{1}{2}+\delta }\right)
\ll B^4\left(M^3N^2q+M^3N^4q^{\frac{1}{2}+\delta}\right)\lnc^{2c_1+c_2}.
\end{align*}
The Lemma now follows from this estimate and (\ref{formula sveden W k W(j,k)}).

{\corollary
\label{Sledst lemmi Dvoynie summi kv sledst}  Let   $M$, $N$, $U$ be integers, $N\le U<2N$,  $q^{\frac{1}{4}-\theta}\le N\le q^{\frac{1}{4}+\theta}$, \linebreak $D^{\frac{\scriptstyle1}{\scriptstyle2}}\le q\le D$, $\nu\le \exp(\sqrt{2\lnc})$,  $a_m$ and  $b_n$ are such integer-valued functions that $|a_m|\le \tau_5(m)$, $|b_n|\le 1$. Then the following estimate holds for $x\ge q^{\frac{3}{4}+\theta+1,1\delta }$:
\begin{align*}
W=\sum_{M<m\le 2M} \hspace{-10pt} a_m\sum_{\substack{U<n\le \min(xm^{-1},2N)\\ (mn,q)=1,\,mn\equiv l\hspace{-8pt} \pmod \nu}} \hspace{-20pt}b_n\chi_q (mn-l)\ll
\frac{x}{\nu}\,\exp\left(-0.7\sqrt{\lnc}\right).
\end{align*}
}

{\sc Proof.} Taking into account that $\ln M\ll \lnc$, we have from the Lemma \ref{Mardjhanashvili}
$$
\sum_{M<m\le 2M}\tau_5(m)\ll M\lnc^4,\qquad \sum_{M<m\le 2M}\tau_5^2(m)\ll M\lnc^{24}.
$$
Applying the Lemma \ref{Lemma dvoynie summi kv osredn} for $c_1=4$, $c_2=24$ and using conditions $MN\le x$,  $q^{\frac{1}{4}-\theta}\le N\le q^{\frac{1}{4}+\theta}$ and $x\ge q^{\frac{3}{4}+\theta+1,1\delta }$, we find that
\begin{align*}
W& \ll \left(M^\frac{3}{4}N^\frac{1}{2}q^{\frac{1}{4}}+M^\frac{3}{4}Nq^{\frac{1}{8}}\right)\lnc^9\delta^\frac{\delta}{4}\le x^\frac{3}{4} \left(N^{-\frac{1}{4}}q^{\frac{1}{4}}+N^{\frac{1}{4}}q^{\frac{1}{8}}\right)\lnc^9q^\frac{\delta}{4}\\
& \ll x\left(\frac{qN^{-1}}{x}+\frac{Nq^{\frac{1}{2}}}{x}\right)^{\frac{1}{4}}\lnc^9q^\frac{\delta}{4}
 \ll x\left(\frac{q^{\frac{3}{4}+\theta}}{x}\right)^{\frac{1}{4}}\lnc^9q^\frac{\delta}{4}\ll x\lnc^9q^{-\frac{\delta}{40}}.
\end{align*}
Applying relations $D^{\frac{\scriptstyle 1}{\scriptstyle 2}}\le q\le D$ and $\nu\le \exp(\sqrt{2\lnc})$, we obtain
$$
W\ll \frac{x}{\nu}\lnc^9\exp(\sqrt{2\lnc})\,D^{-\frac{\delta}{80}}\ll \frac{x}{\nu}\,\exp\left(-0.7\sqrt{\lnc}\right).
$$

{\lemma \label{Lemma dvoynie summi cub osredn} Let  $M$, $N$, $U$ be integers, $N\le U<2N\le  q^{\frac{1}{6}}$, $a_m$  and $b_n$ are such integer-valued functions that
$$
\sum_{M<m\le 2M}|a_m|^\alpha\ll M\lnc^{c_\alpha},\quad  \alpha =1,2; \qquad |b_n|\ll B.
$$
Then the following estimate holds
\begin{align*}
W=\sum_{M<m\le 2M} \hspace{-10pt} a_m\sum_{\substack{U<n\le \min(xm^{-1},2N)\\ (mn,q)=1,\,mn\equiv l\hspace{-8pt} \pmod \nu}} \hspace{-20pt}b_n\chi_q (mn-l)\ll  BM^{\frac 56}N^{\frac{1}{2}}q^{\frac{1}{6}+\frac{1}{6}\delta }\lnc^{\frac{4c_1+c_2}{6}+1}.
\end{align*}}

\textsc{Proof.}  Without loss of generality, we can assume that $MN<x$. Repeating the argument used in the proof of the previous lemma, we find
\begin{align}
|W|&\ll \lnc \max_{0\le j<\nu}\max_{0\le k<q}W(j,k), \qquad W(j,k)=\sum_{\substack{M<m\le 2M\\ (m,q)=1}}|a_m||\B(k\nu+jq,m)|. \label{formula sveden W k W(j,k)-3}
\end{align}
\begin{align*}
\B(k\nu+jq,m)=\sum_{\substack{N<n\le 2N\\ (n,q)=1}}b_n\chi_q(mn-l)e\left(\frac{(k\nu+jq)n}{q\nu}\right).
\end{align*}
Let us estimate $W(j,k)$. Cubing both sides of the identity and applying H\"{o}lder inequality, we have
\begin{align*}
W^3(j,k)&\ll M^2\lnc^{2c_1} \sum_{\substack{M<m\le 2M\\ (m,q)=1}}\hspace{-6pt}|a_m||\B(kd+jq,m)|^3.
\end{align*}
Now, squaring both sides of the last estimate and applying Cauchy inequality, we find
\begin{align*}
W^6(j,k)&\ll M^4\lnc^{4c_1}\sum_{\substack{M<m\le 2M\\ (m,q)=1}}|a_m|^2\sum_{\substack{M<m\le 2M\\ (m,q)=1}}|\B(kd+jq,m)|^6\\
&\ll M^5\lnc^{4c_1+c_2}\sum_{\substack{m=0\\(m,q)=1}}^{q-1}\left|\sum_{\substack{N<n\le 2N\\ (n,q)=1}}b_n\chi_q(n-lm_q^{-1})e\left(\frac{(kd+jq)n}{qd}\right)\right|^6\\
&=M^5\lnc^{4c_1+c_2}\sum_{\substack{\lambda=0\\ (\lambda,q)=1}}^{q-1}\left|\sum_{\substack{N<n\le 2N\\ (n,q)=1}}\hspace{-6pt}b_n\chi_q(n+\lambda)e\left(\frac{(kd+jq)n}{qd}\right)\right|^6\\
&\ll M^5\lnc^{4c_1+c_2}\hspace{-12pt}\sum_{\substack{N'<n_1,\ldots,n_6\le 2N\\ (n_1,\ldots,n_6),q)=1}}\hspace{-8pt}|b_{n_1}\ldots b_{n_6}|
\left|\sum_{\lambda =0}^{q-1}\chi \left(\frac{(\lambda +n_1)(\lambda +n_2)(\lambda +n_3)}{(\lambda +n_4)(\lambda +n_5)(\lambda +n_6)}\right)\right|.
\end{align*}
Next, using the estimate $|b_n|\ll B$ and the Lemma \ref{Teorema Burgess cub}, we find that
\begin{align*}
W^6(j,k) &
\ll B^6M^5N^3q^{1+\delta}\lnc^{4c_1+c_2}.
\end{align*}
The Lemma now follows from this estimate and (\ref{formula sveden W k W(j,k)-3}).

{\corollary
 \label{Sledst lemmi Dvoynie summi cub sledst}  Let $M$, $N$, $U$ be integers, $N\le U<2N$,  $q^{\theta}\le N\le q^{\frac{1}{6}}$, $D^{\frac{\scriptstyle 1}{\scriptstyle 2}}\le q\le D$, $\nu\le \exp(\sqrt{2\lnc})$,  $a_m$ and  $b_n$ are such integer-valued functions that $|a_m|\le \tau_5(m)$, $|b_n|\le 1$. Then the following estimate holds for $x\ge q^{1-2\theta+1,1\delta }$:
\begin{align*}
W=\sum_{M<m\le 2M} \hspace{-10pt} a_m\sum_{\substack{U<n\le \min(xm^{-1},2N)\\ (mn,q)=1,\,mn\equiv l\hspace{-8pt} \pmod \nu}} \hspace{-20pt}b_n\chi_q (mn-l)\ll
\frac{x}{\nu}\,\exp\left(-0.7\sqrt{\lnc}\right).
\end{align*}
}

{\sc Proof.}  Taking into account that $\ln M\ll \lnc$, we have from the Lemma \ref{Mardjhanashvili} that
$$
\sum_{M<m\le 2M}\tau_5(m)\ll M\lnc^4,\qquad \sum_{M<m\le 2M}\tau_5^2(m)\ll M\lnc^{24}.
$$
Applying the Lemma $x~\ge~q^{1-2\theta +1,1\delta }$ for $c_1=4$, $c_2=24$  and using conditions $MN\le x$, $N\ge q^\theta$ and $x~\ge~q^{1-2\theta +1,1\delta }$, we find that
\begin{align*}
W&\ll (MN)^{\frac 56}N^{-\frac{1}{3}}q^{\frac{1}{6}+\frac{1}{6}\delta }\lnc^\frac{20}{3}\le  x^{\frac 56}N^{-\frac{1}{3}}q^{\frac{1}{6}+\frac{1}{6}\delta }\lnc^\frac{20}{3}= \\
&=x  \left(\frac{N^{-2}q^{1+\delta }}x\right)^{\frac{1}{6}}\lnc^\frac{20}{3}\le x\left(\frac{q^{1-2\theta+\delta }}x\right)^{\frac{1}{6}}\lnc^\frac{20}{3}\ll xq^{-\frac{\delta}{60}}\lnc^\frac{20}{3}.
\end{align*}
Applying relations $D^{\frac{\scriptstyle 1}{\scriptstyle 2}}\le q\le D$ and $\nu\le \exp(\sqrt{2\lnc})$, we obtain
$$
W\ll \frac{x}{\nu}\lnc^\frac{20}{3}\exp(\sqrt{2\lnc})\,D^{-\frac{\delta}{120}}\ll \frac{x}{\nu}\,\exp\left(-0.7\sqrt{\lnc}\right).
$$


\section{Proof of the Theorem \ref{Teorema ocnovnaya}}
Without loss of generality, we shall assume that
$$
x=D^{\frac{\scriptstyle 5}{\scriptstyle 6}+\varepsilon}\ge q^{\frac{\scriptstyle 5}{\scriptstyle 6}+\varepsilon}.
$$
Having in mind that the contribution of terms satisfying $(n,q)>1$ in the sum $T(\chi)$ has order of magnitude $\ll \lnc^2$ and using that $\chi_q$ is a primitive character generated by a non-principal character $\chi$ and $q_1$ is a product of primes, dividing $D$, but not $q$, we have
\begin{align*}
T(\chi)&=\sum_{\substack{n\le x\\ (n,q)=1}}\Lambda (n)\chi (n-l)+O(\lnc^2)=\sum_{\substack{n\le x,\ (n,q)=1\\ (n-l,q_1)=1}}\Lambda (n)\chi_q(n-l)+O(\lnc^2)\\ &=\sum_{\nu| q_1}\mu (\nu)T(\chi_q,\nu)+O(\lnc^2),\qquad T(\chi_q,\nu)=\sum_{\substack{n\le x,\ (n,q)=1\\ n\equiv l\hspace{-8pt}\pmod \nu}}\Lambda(n)\chi_q(n-l).
\end{align*}
Part of the sum $T(\chi )$ corresponding to terms satisfying the condition $\exp(\sqrt{2\lnc})<\nu\le x$ shall be denoted by $T_1(\chi )$. Let us estimate $T_1(\chi )$ using the trivial estimate for the sum $T(\chi_q,\nu)$ and the Lemma \ref{Lemma summs s bolshimi delit}:
\begin{align*}
&|T_1(\chi )|\ll \lnc\hspace{-8pt}\sum_{\substack{\nu| q_1\\ \nu>\exp \sqrt{2\lnc}}}\hspace{-8pt}\mu^2(\nu)\left(\frac{x}{\nu}+1\right)\ll x\lnc\hspace{-8pt}\sum_{\substack{\nu| q_1\\ \nu>\exp \sqrt{2\lnc}}}\hspace{-8pt}\frac{\mu^2 (\nu)}{\nu}\ll x\lnc \exp\left(-0.7\sqrt{\lnc}\right).
\end{align*}
Therefore
\begin{align}\label{formula T(chi)=sumT(chi,nu)}
 T(\chi )=&\sum_{\substack{\nu|  q_1\\ \nu\le\exp(\sqrt{2\lnc})}}\hspace{-14pt}\mu(\nu)T(\chi_q,\nu)+O\left(x\lnc\exp\left(-0.7\sqrt{\lnc}\right)\right).
\end{align}
Let us estimate $T(\chi_q,\nu)$ for $\nu\le \exp \sqrt{2\lnc}$ and $(\nu,l)=(q,l)=(\nu,q)=1$. The sum $T(\chi_q,\nu)$ is studied in detail in the Lemma 5 of \cite{RakhmonovZKh-1994-TrMIRAN} and we use the following estimate
$$
|T(\chi_q,\nu)|\le 10x \ln^5x\left( \sqrt{\frac{1}{q\nu^2}+\frac{q}{x}} +x^{-\frac{1}{6}}\nu^{-\frac{1}{2}}+ x^{-\frac{1}{3}}q^\frac{1}{6}\nu^{-\frac{1}{3}}\right)\tau(q).
$$
For $x>q^{1+1,2\varepsilon}$ and $q\ge \exp \sqrt{2\lnc}$ the last estimate gives the following non-trivial estimate for the sum
$$
|T(\chi_q,\nu)|\ll \frac{x}{\nu}\exp\left(-0.7\sqrt{\lnc}\right).
$$
Thus, henceforth we shall assume that $x\le q^{1+1,2\varepsilon}$, i.e.
\begin{equation}\label{formulaD^{5/6}<q}
D^{\frac{\scriptstyle5}{\scriptstyle6}}\le q\le D.
\end{equation}
Setting $u=x^{\frac{1}{3}}$, $r=3$ in the Lemma \ref{Resheto} and
$$
f(n)=\left\{
       \begin{array}{ll}
        \chi_q (n-l), & \hbox{ for $(n,q)=1$ and $n\equiv l\hspace{-8pt}\pmod \nu$;} \\
         0. & \hbox{ otherwise,}
       \end{array}
     \right.
$$
we find
\begin{equation}
\label{formula T(chi,nu)=sum_Tk(chi,nu)}
T(\chi_q,\nu)=\sum_{k=1}^3(-1)^kC_3^k\tilde{T}_k(\chi_q,\nu),
\end{equation}
$$
\tilde{T}_k(\chi_q,\nu)=\sum_{m_1\le u}\mu (m_1)\cdots \sum_{m_k\le u}\mu (m_k)\hspace{-90pt}\sum_{\substack {n_1 \\ m_1\cdots m_kn_1\cdots n_k\le x, \ (m_1\cdots m_kn_1\cdots n_k,q)=1,\ m_1\cdots m_kn_1\cdots n_k\equiv l\hspace{-8pt} \pmod \nu}} \hspace{-90pt}\ldots \sum_{n_k}\ln n_1\chi_q(m_1n_1\cdots m_kn_k-l),
$$
Let us divide in $T_k(\chi_q,\nu)$ the limits of each variable $m_1,\cdots ,m_k,n_1, \cdots ,n_k$ into not more than $\lnc$ intervals of the form $M_j<m_j\le 2M_j$, $N_j<n_j\le 2N_j$, $j=1,2,\cdots ,k$. We obtain not more than $\lnc^{2k}$  sums of the form
\begin{align*}
\hat{T}_k&(\chi_q,\nu)=\hspace{-10pt}\sum_{M_1<m_1\le 2M_1}\hspace{-15pt}\mu (m_1)\cdots\hspace{-16pt}\sum_{M_k<m_k\le 2M_k}\hspace{-15pt}\mu (m_k) \hspace{-80pt}\sum_{\substack{N_1<n_1\le 2N_1\\ m_1n_1\cdots m_kn_k\le x, \ (m_1n_1\cdots m_kn_k,q)=1,\ m_1n_1\cdots m_kn_k\equiv l\hspace{-8pt} \pmod \nu}}\hspace{-80pt}\cdots\hspace{-2pt}\sum_{N_k<n_k\le 2N_k}\hspace{-3mm} \hspace{-4pt}\chi_q (m_1n_1\cdots m_kn_k-l)\ln n_1\\
&=\int\limits_1^{2N_1}\hspace{-5pt}\sum_{M_1<m_1\le 2M_1}\hspace{-15pt}\mu (m_1)\cdots\hspace{-16pt}\sum_{M_k<m_k\le 2M_k}\hspace{-15pt}\mu (m_k) \hspace{-70pt}\sum_{\substack{\max (u,N_1)<n_1\le 2N_1\\ m_1n_1\cdots m_kn_k\le x, \ (m_1n_1\cdots m_kn_k,q)=1,\ m_1n_1\cdots m_kn_k\equiv l\hspace{-8pt} \pmod \nu}}\hspace{-70pt}\cdots\hspace{-2pt}\sum_{N_k<n_k\le 2N_k}\hspace{-3mm} \hspace{-4pt}\chi_q (m_1n_1\cdots m_kn_k-l)d\ln u.
\end{align*}
Let us denote by $U_1=\max (u,N_1)$ such number $u$ for which the integrand takes on its maximum value. Then we have the following inequality $|\hat{T}_k(\chi_q,\nu)|\ll \lnc \left|T_k(\chi_q,\nu)\right|$, where
\begin{align}
T_k&(\chi_q ,\nu)=\hspace{-10pt}\sum_{M_1<m_1\le 2M_1}\hspace{-15pt}\mu(m_1)\cdots\hspace{-80pt}\sum_{\substack {M_k<m_k\le 2M_k \\ m_1n_1\cdots m_kn_k\le x, \ (m_1n_1\cdots m_kn_k,q)=1,\ m_1n_1\cdots m_kn_k\equiv l\hspace{-8pt} \pmod \nu}}\hspace{-80pt}\mu (m_k)\hspace{-5pt}\sum_{U_1<n_1\le 2N_1}\hspace{-5pt}\cdots \hspace{-5pt}\sum_{U_k<n_k\le 2N_k}\hspace{-10pt}\chi_q (m_1n_1\cdots m_kn_k-l), \nonumber \\
&x^{\frac{1}{3}}>M_1\ge M_2\ge\cdots \ge M_k, \qquad \hspace{-10pt}N_1\ge N_2\ge\cdots \ge N_k,\qquad \hspace{-20pt} N_j\le U_j<2N_j.\label{usl M_1 ge M_2 ge ... ge M_k, N_1 ge N_2 ge ... ge N_k}
  \end{align}
  From this and the estimates (\ref{formula T(chi,nu)=sum_Tk(chi,nu)}), (\ref{formula T(chi)=sumT(chi,nu)}), we obtian
\begin{equation}\label{formula T(chi)=summTk(chi,q,d)}
|T(\chi)|\ll\sum_{\substack{\nu| q_1\\\nu\le\exp(\sqrt{2\lnc})}}\hspace{-10pt}\mu^2(\nu)\sum_{k=1}^3\lnc^{6}\max|T_k(\chi_q,\nu)| +x\exp\left(-0.6\sqrt{\lnc}\right).
\end{equation}
Let us introduce the following notations:
\begin{align*}
  & \prod_{j=1}^kM_jN_j=Y, \qquad \prod_{j=1}^kM_jU_j=X, \qquad Y<X\le x,
  \end{align*}
and further we shall assume that
\begin{equation}
\label{formula usl Y>x exp(...)}
Y\ge x\exp\left(-\sqrt{\lnc}\right),
\end{equation}
since otherwise, estimating $T_k(\chi_q,\nu)$ trivially, we have
\begin{align*}
T_k(\chi_q,\nu)&\ll \sum_{\substack{X<n\le 2^kY\\n\equiv l\hspace{-8pt} \pmod \nu}}\tau_{2k}(n)\ll\frac{ 2^kY}{\nu}\lnc^{2k-1}\\
&\le\frac{x}{\nu}2^k\lnc^{2k-1}\exp\left(-\sqrt{\lnc}\right)\ll \frac{x}{\nu}\exp\left(-0.7\sqrt{\lnc}\right).
\end{align*}
The sums $T_k(\chi_q,\nu)$, $k=1,2,3$ are estimated in the same way. For example, we shall estimate the sum $T_3(\chi_q,\nu)$ and consider the following possible values for the parameter $N_1$:
\begin{enumerate}
\item  $N_1>q^{\frac{\scriptstyle1}{\scriptstyle3}+\frac{\scriptstyle{16}}{\scriptstyle{27}}\varepsilon}$;
\item   $q^\frac{\scriptstyle1}{\scriptstyle6}< N_1\le q^{\frac{\scriptstyle1}{\scriptstyle3}+\frac{\scriptstyle{16}}{\scriptstyle{27}}\varepsilon}$;
\item  $q^\frac{\scriptstyle1}{\scriptstyle{12}}<N_1\le q^\frac{\scriptstyle1}{\scriptstyle6}$,
\item   $N_1<q^\frac{\scriptstyle1}{\scriptstyle{12}}$.
\end{enumerate}
In considering the cases 1, 2, we shall transform the sum $T_3(\chi_q,\nu)$ and rewrite it in the form
\begin{align*}
T_3(\chi_q,\nu)&=\sum_{XU_1^{-1}<m\le 2^{5}YN_1^{-1}}\hspace{-20pt} a_m \sum_{\substack{U_1<n\le 2N_1,\,mn\le x\\(mn,q)=1,\, mn\equiv l\hspace{-8pt}\pmod \nu}}\hspace{-20pt} \chi_q (mn-l), \qquad |a_m|\le \tau_{5}(m), 
\end{align*}
Next, we shall divide the interval of summation $XU_1^{-1}<m\le 2^{5}YN_1^{-1}$ into smaller intervals of the form  $M<m\le 2M$. We shall obtain not more than $5$ sums of the form
$$
T_3(\chi_q,\nu,M)=\sum_{M<m\le 2M} \hspace{-10pt} a_m\sum_{\substack{U_1<n\le \min(xm^{-1},2N_1)\\ (mn,q)=1,\,mn\equiv l\hspace{-8pt} \pmod \nu}} \hspace{-10pt}\chi_q (mn-l).
$$

\textbf{Case 1.} $N_1>q^{\frac{\scriptstyle1}{\scriptstyle3}+\frac{\scriptstyle{16}}{\scriptstyle{27}}\varepsilon}$. Determining the number $m_q^{-1}$ from the congruence $mm_q^{-1}\equiv 1 \pmod q$ and passing to the estimate, we find
\begin{align*}
|T_3(\chi_q,\nu,M)|&\le \sum_{\substack{M<m\le 2M \\ (m,q)=1}}\tau_{5}(m)\left|\sum_{\substack{U_1<n\le\min(xm^{-1},2N_1)\\ (n,q)=1,\,mn\equiv l\hspace{-8pt} \pmod \nu}}\chi_q(mn-l)\right|.
\end{align*}
Applying the Lemma \ref{Lemma korotkaya summa} to the sum over $n$ for
$$
\delta=\frac{10}{27}\varepsilon,\quad \eta=lm_q^{-1},\quad u=\min(xm^{-1},2N_1),\quad y=\min(xm^{-1},2N_1)-U_1\le N_1,
$$
and using the condition $N_1>q^{\frac{\scriptstyle1}{\scriptstyle3}+\frac{\scriptstyle{16}}{\scriptstyle{27}}\varepsilon}$, we have
\begin{align*}
|T_3(\chi_q,\nu,M)|&\ll M\lnc^4\cdot\frac{N_1}{\nu}\exp\left(-0.7\sqrt{\lnc}\right)\\
&\le\frac{2^{5}Y}{\nu}\lnc^4\exp\left(-0.7\sqrt{\lnc}\right)\ll\hspace{-5pt}\ \frac{x}{\nu}\lnc^4\exp\left(-0.7\sqrt{\lnc}\right).
\end{align*}

\textbf{Case 2. $q^{\frac{\scriptstyle1}{\scriptstyle6}}<N_1\le q^{\frac{\scriptstyle1}{\scriptstyle3}+\frac{\scriptstyle{16}}{\scriptstyle{27}}\varepsilon}$.} We shall use the Corollary \ref{Sledst lemmi Dvoynie summi kv sledst} of the Lemma \ref{Lemma dvoynie summi kv osredn} for
$$
U=U_1, \qquad N=N_1, \qquad \theta =\frac{1}{12}+\frac{16}{27}\varepsilon,\qquad \delta=\frac{10}{27}\varepsilon, \qquad b_n=1.
$$
Then, for $x\ge q^{\frac{\scriptstyle3}{\scriptstyle4}+\theta +1,1\delta}= q^{\frac{\scriptstyle5}{\scriptstyle6}+\varepsilon}$ we have
\begin{align*}
|T_3(\chi_q,\nu,M)|\ll  \frac{x}{\nu}\,\exp\left(-0.7\sqrt{\lnc}\right).
\end{align*}

\textbf{Case 3. $q^\frac{\scriptstyle1}{\scriptstyle{12}}<N_1\le q^{\frac{\scriptstyle1}{\scriptstyle6}}$}.The conditions of the Corollary \ref{Sledst lemmi Dvoynie summi cub sledst} are satisfied for the sum $T_3(\chi_q,\nu,M)$  for
$$
U=U_1,\qquad N=N_1,  \qquad \theta =\frac{1}{12}, \qquad b_n=1, \qquad \delta=\frac{10}{27}\varepsilon.
$$
Applying this Corollary for $x\ge q^{\scriptstyle1-\scriptstyle{2\theta}+\scriptstyle{1,1\delta}}= q^{\frac{5}{6}+\frac{\scriptstyle{11}}{\scriptstyle{27}}\varepsilon}$, we obtain
\begin{align*}
|T_3(\chi_q,\nu,M)|\ll  \frac{x}{\nu}\,\exp\left(-0.7\sqrt{\lnc}\right).
\end{align*}

\textbf{Case 4.  $N_1<q^\frac{\scriptstyle1}{\scriptstyle{12}}$}. We shall transform the sum $T_3(\chi_q,\nu)$ by rewriting it in the form
\begin{align*}
T_3(\chi_q,\nu)&=\sum_{XM_1^{-1}<m\le 2^{5}YM_1^{-1}}\hspace{-10pt} a_m \sum_{\substack{M_1<n\le 2M_1,\,mn\le x\\(mn,q)=1,\, mn\equiv l\hspace{-3pt}\pmod \nu}}\hspace{-20pt}\mu(n)\chi_q (mn-l), \qquad |a_m|\le \tau_{5}(m),
\end{align*}
and divide the interval of summation $XM_1^{-1}<m\le 2^{5}YM_1^{-1}$ into smaller intervals of the form  $M<m\le 2M$. We shall obtain not more than $5$ sums of the form
$$
T_3(\chi_q,\nu,M)=\sum_{M<m\le 2M} \hspace{-10pt} a_m\sum_{\substack{M_1<n\le \min(xm^{-1},2M_1)\\ (mn,q)=1,\,mn\equiv l\hspace{-8pt} \pmod \nu}} \hspace{-10pt}\mu(n)\chi_q (mn-l).
$$
Using the relations (\ref{usl M_1 ge M_2 ge ... ge M_k, N_1 ge N_2 ge ... ge N_k}), (\ref{formula usl Y>x exp(...)}),  conditions of the case 4 and the relations (\ref{formulaD^{5/6}<q}),  we have
\begin{align*}
M_1&\ge \left(M_1M_2M_3\right)^{\frac{\scriptstyle1}{\scriptstyle3}}= \left(\frac{Y}{N_1N_2N_3}\right)^{\frac{\scriptstyle1}{\scriptstyle3}}\ge \frac{Y^{\frac{\scriptstyle1}{\scriptstyle3}}}{N_1}
\ge \frac{\left(x\,\exp\left(-\sqrt{\lnc}\right)\right)^{\frac{\scriptstyle1}{\scriptstyle3}}}{N_1}\\
& \qquad \qquad \ge\frac{q^{\frac{\scriptstyle5}{\scriptstyle{18}}+\frac{\scriptstyle5}{\scriptstyle{18}}\varepsilon}\, \exp\left(-0.4\sqrt{\lnc}\right)}{q^{\frac{\scriptstyle1}{\scriptstyle{12}}}}>q^{\frac{\scriptstyle7}{\scriptstyle{36}}},\\
M_1&\le x^{\frac{1}{3}}=D^{\frac{\scriptstyle5}{\scriptstyle{18}}+\frac{\scriptstyle\varepsilon}{\scriptstyle3}}\le q^{\frac{\scriptstyle1}{\scriptstyle3}+\frac{\scriptstyle2}{\scriptstyle5}\scriptstyle\varepsilon}.
\end{align*}
It follows that the conditions of the Corollary \ref{Sledst lemmi Dvoynie summi kv sledst} are satisfied for the sum $T_3(\chi_q,\nu,M)$ for
$$
U=M_1, \qquad N=M_1,\qquad \theta =\frac{1}{12}+\frac{2}{5}\varepsilon, \qquad \delta=\frac{10}{27}\varepsilon, \qquad b_n=\mu(n)
$$
Applying this Corollary for
$$
x\ge q^{\frac{\scriptstyle3}{\scriptstyle4}+\theta +1,1\delta}= q^{\frac{\scriptstyle5}{\scriptstyle6}+\frac{\scriptstyle{109}}{\scriptstyle{135}}\varepsilon},
$$
we obtain
\begin{align*}
|T_3(\chi_q,\nu,M)|\ll  \frac{x}{\nu}\,\exp\left(-0.7\sqrt{\lnc}\right).
\end{align*}
The Lemma now follows by inserting the estimates $T_k(\chi_q,\nu)$, $k=1,2,3$ into (\ref{formula T(chi)=summTk(chi,q,d)}).


\vspace{10mm}
\begin{thebibliography}{12} 

\bibitem{Burgess-Proc.Lon.M.S (3)12-1962} D.~A.~Burgess, \emph{On character sums and $L$~--~series}, Proc. London Math. Soc. 1962, v. 12, no 3, pp. 193-206.
\bibitem{Burgess-1986}  D.~A.~Burgess, \emph{The character sum estimate with $r = 3$}, J. London Math. Soc. \textbf{33} (1986). 219-226.
\bibitem{Friedlander-Gong-Shparlinski} J.~B.~Friedlander, K.~Gong, I.~E.~Shparlinski, \emph{Character Sums over Shifted Primes}, Math. Notes, \textbf{88:3-4} (2010), 585-598. 
\bibitem{Jutila}  M.~Jutila, \emph{On the least Goldbach's number in an arithmetical progression with a prime difference}, Ann. Univ. Turku; Ser. A., I, \textbf{118} (1968).
\bibitem{Heath-Brown-1982} D.~R.~Heath-Brown, \emph{Prime numbers in short intervals and a generalized Vaughan identity}, Canad. J. Math. 34, 1982, 1365-1377.
\bibitem{Huxley-1942-Inv.Math} M.~N.~Huxley,  \emph{On the difference between consecutive primes}, Inventiones mathematical, June, \textbf{15:2}  (1971), 164-170.
\bibitem{KaratsubaAA-2008} A.~A.~Karatsuba, \emph{Arithmetic problems in the theory of Dirichlet characters}, Russian Math. Surveys, \textbf{63:4} (2008), 641-690.
\bibitem{KaratsubaAA-1968} A.~A.~Karatsuba, \emph{Sums of characters, and primitive roots, in finite fields}, Dokl. AN SSSR, \textbf{180} (1968), 1287-1289, (Russian).
\bibitem{KaratsubaAA-1970-1} A.~A.~Karatsuba, \emph{Estimates of character sums}, Math. USSR-Izv., \textbf{4:1} (1970), 19-29.
\bibitem{KaratsubaAA-1970-2} A.~A.~Karatsuba, \emph{Sums of characters over prime numbers}, Math. USSR-Izv., \textbf{4:2} (1970), 303-326.
\bibitem{KaratsubaAA-1970-3} A.~A.~Karatsuba, \emph{Sums of characters with prime numbers}, Dokl. AN SSSR, \textbf{190} (1970), 517-518, (Russian).
\bibitem{KaratsubaAA-1970-4} A.~A.~Karatsuba, \emph{The distribution of products of shifted prime numbers in arithmetic progressions}, Dokl. AN SSSR, \textbf{192} (1970), 724-727, (Russian).
\bibitem{KaratsubaAA-1971} A.~A.~Karatsuba, \emph{ Sums of characters with prime numbers in an arithmetic progression}, Math. USSR-Izv., \textbf{5:3} (1971), 485-501.
\bibitem{KaratsubaAA-1975-1} A.~A.~Karatsuba, \emph{Sums of characters in sequences of shifted prime numbers, with applications}, Math. Notes, \textbf{17:1} (1975), 91-93.
\bibitem{KaratsubaAA-1975-2} A.~A.~Karatsuba, \emph{Some problems of contemporary analytic number theorem}, Math. Notes, \textbf{17:2} (1975), 195-199 .
\bibitem{KaratsubaAA-1976} A.~A.~Karatsuba,  \emph{Distribution of values of nonprincipal characters}, Proc. Steklov Inst. Math., \textbf{142} (1979), 165-174.
\bibitem{KaratsubaAA-1978} A.~A.~Karatsuba, \emph{Sums of Legendre symbols of polynomials of second degree over prime numbers}, Math. USSR-Izv., \textbf{12:2} (1978), 299-308.
\bibitem{KaratsubaOATCh}A.~A.~Karatsuba, \emph{Osnovy analiticheskoi teorii chisel}, 2-e izd., Nauka, M., 1983 , 240 pp., (Russian).
\bibitem{Kerr-2017} B.~Kerr, \emph{On certain exponential and character sums}, PhD Thesis, UNSW, 2017.
\bibitem{Linnik-1952} Ju.~V.~Linnik, \emph{Some conditional theorems conceming binary problems with prime numbers,} Izv. Akad. Nauk SSSR Ser. Mat., \textbf{16} (1952), 503-520. (Russian).
\bibitem{Linnik-1973} Ju.~V.~Linnik, \emph{Recent works of I.~M.~Vinogradov}, Trudy Mat. Inst. Steklov., \textbf{132}, 1973, 27-29; Proc. Steklov Inst. Math., \textbf{132} (1975), 25-28.
\bibitem{Mardjhanashvili} K.~K.~Mardjhanashvili, \emph{An estimate for an arithmetic sum},  Doklady Akad. Nauk SSSR, \textbf{22:7}, 391-393, (Russian).
\bibitem{Prachar_Acta-Arith-29(1976)} K.~Prachar, \emph{Uber die Anwendung einer Methode von Linnik}, Acta Arith. \textbf{29} (1976), 367-376.
\bibitem{Prachar_Acta-Arith-44(1984)} K.~Prachar, \emph{Bemerkungen uber Primzahlen in kurzen Reihen. [Remarks on primes in short sequences]}, Acta Arith. \textbf{44} (1984), 175-180.
\bibitem{RakhmonovZKh-1986-UMN} Z.~Kh.~Rakhmonov \emph{On the distribution of values of Dirichlet characters},
Russian Math. Surveys 41 (1986), no. 1, 237-238.  
\bibitem{RakhmonovZKh-1986-DANRT} Z.~Kh.~Rakhmonov, \emph{Estimation of the sum of characters with primes}, Dokl. Akad. Nauk Tadzhik. SSR, \textbf{29:1}, 16-20. (Russian).
\bibitem{RakhmonovZKh-1994-TrMIRAN} Z.~Kh.~Rakhmonov, \emph{On the distribution of the values of Dirichlet characters and their applications}, 
    Proc. Steklov Inst. Math.,  \textbf{207:6} (1995),  263-272.
\bibitem{RakhmonovZKh-1986-IzvANRT} Z.~Kh.~Rakhmonov, \emph{The least Goldbach number in an arithmetic progression}, Izv. Akad. Nauk Tadzhik. SSR Otdel. Fiz.-Mat. Khim. i Geol. Nauk  \textbf{2(100)}, 1986, 103-106, (Russian).
\bibitem{RakhmonovZKh-2013-DANRT} Z.~Kh.~Rakhmonov, \emph{Distribution of values of Dirichlet characters in the sequence of shifted primes}, Doklady Akademii nauk Respubliki Tajikistan,  \textbf{56:1} (2013), 5-9, (Russian).
\bibitem{RakhmonovZKh-2013-IzvSarUniv} Z.~Kh.~Rakhmonov, \emph{Distribution of values of Dirichlet characters in the sequence of shifted primes}, Izv. Saratov Univ. (N.S.), Ser. Math. Mech. Inform., 13:4(2) (2013),  113-117, (Russian).
\bibitem{RakhmonovZKh-2014-Ch.sbor-15-2(50)} Z.~Kh.~Rakhmonov, \emph{	Sums of characters over prime numbers}, Chebyshevskii Sb., \textbf{15:2} (2014), 73-100. (Russian).
\bibitem{RakhmonovZKh-2017-TrMIRAN-299} Z.~Kh.~Rakhmonov, \emph{The sums of the values of nonprincipal characters from the sequence of shifted prime}, Trudy Mat. Inst. Steklov, \textbf{299} (2017), 1-27, (Russian).
\bibitem{RakhmonovZKh-2017-DANRT-9} Z.~Kh.~Rakhmonov, \emph{On the estimation of the sum the values of Dirichlet character in a sequence of shifted primes}, Doklady Akademii nauk Respubliki Tajikistan,  \textbf{60:9} (2017), 378-382, (Russian).
\bibitem{RakhmonovZKh-1993-IzvRAN} Z.~Kh.~Rakhmonov \emph{A theorem on the mean value of $\psi (x,\chi)$ and its applications},
\bibitem{VinogradovAI} A.~I.~Vinogradov, \emph{On numbers with small prime divisors}, Dokl. Akad. Nauk SSSR, \textbf{109:4} (1956), 683-686, (Russian).
\bibitem{VinigradovIM-1938} I.~M.~Vinogradow, \textit{Some general lemmas and their application to the estimation of trigonometrical sums}, Rec. Math. [Mat. Sbornik] N.S., \textbf{3(45):3} (1938), 435-471, (Russian).
\bibitem{VinigradovIM-1943} I.~M.~Vinogradow, \textit{An improvement of the estimation of sums with primes}, Izv. Akad. Nauk SSSR Ser. Mat., \textbf{7:1} (1943), 17-34,(Russian).
\bibitem{VinigradovIM-1952} I.~M.~Vinogradow, \textit{New approach to the estimation of a sum of values of $\chi(p+k$)}, Izv. Akad. Nauk SSSR Ser. Mat., \textbf{16:3} (1952), 197-210. (Russian).
\bibitem{VinigradovIM-1953}I.~M.~Vinogradow, \textit{Improvement of an estimate for the sum of the values $\chi (p+k)$}, Izv. Akad. Nauk SSSR Ser. Mat., \textbf{17:4} (1953), 285-290. (Russian).
\bibitem{VinigradovIM-1966}I.~M.~Vinogradow, \textit{An estimate for a certain sum extended over the primes of an arithmetic progression}, Izv. Akad. Nauk SSSR Ser. Mat., \textbf{30:3} (1966), 481-496, (Russian).
\bibitem{WangYuan_SciSinica-1977-20} Wang~Yuan, \emph{On Linnik's method concerning the Goldbach number},   Sci. Sinica, \textbf{20} (1977), 16-30.
\end{thebibliography}
\end{document}